\documentclass[final]{siamltex}

\usepackage{amssymb,latexsym,amsmath}
\usepackage{epsfig}
\usepackage{caption}
 
 \usepackage{epsfig}

\usepackage{makeidx}

\allowdisplaybreaks

 \usepackage{graphicx,fancybox,latexsym,epsfig}
\usepackage{fancyhdr,amsmath,times,amsxtra,amssymb}
\usepackage{color}

\usepackage{fancybox}
\usepackage{tikz}

\usepackage{amsmath}
\usepackage{amsfonts}
\usepackage{array}
\usepackage{dsfont}
\usepackage{hyperref}
\usepackage{amssymb}
\usepackage{bbold}
\usepackage{standalone}
\usepackage{accents}
\usepackage{enumitem}
\usepackage{tikz}
\usepackage{pgfplots}
\pgfplotsset{compat=1.6}
\usepgfplotslibrary{groupplots}
\usepackage{cancel}
 \usepackage{color}

\allowdisplaybreaks

\newtheorem{prop}{Proposition}[section]
\newtheorem{assumption}{Assumption}[section]

  \newtheorem{remark}{Remark}
  \newtheorem{example}{Example}
   \newtheorem{condition1}{Condition (C)}

\usepackage{amssymb,latexsym,amsmath}
\usepackage{mathrsfs}
\usepackage{epsfig}
\usepackage{caption}

\newcommand\scalemath[2]{\scalebox{#1}{\mbox{\ensuremath{\displaystyle #2}}}}

\newcommand{\qedwhite}{\hfill \ensuremath{\Box}}

\pgfplotsset{soldot/.style={color=blue,only marks,mark=*}}
\pgfplotsset{holdot/.style={color=blue,fill=white,only marks,mark=*}}
\fussy

\allowdisplaybreaks

\begin{document}
\sloppy
\title{Isomorphism Properties of Optimality and Equilibrium Solutions under Equivalent Information Structure Transformations I: Stochastic Dynamic Teams
\thanks{Research of first and third authors was supported in part by
the Natural Sciences and Engineering Research Council (NSERC) of Canada. Research of second author was supported in part by an AFOSR Grant (FA9550-19-1-0353). 
	Sina Sanjari and Serdar Y\"uksel are with the Department of Mathematics and 	Statistics, Queen's University, Kingston, ON, Canada,
     Email: \{s.sanjari,yuksel@queensu.ca\},
     Tamer Ba\c{s}ar is with the Coordinated Science Laboratory, University of 		Illinois at Urbana-Champaign, Urbana, IL 61801 USA. E-mail:\{basar1@illinois.edu\}}}
\author{Sina~Sanjari, Tamer Ba\c{s}ar, and Serdar Y\"uksel}
\maketitle

\begin{abstract}
In stochastic optimal control, change of measure arguments have been crucial for stochastic analysis. Such an approach is often called static reduction in dynamic team theory (or decentralized stochastic control) and has been an effective method for establishing existence and approximation results for optimal policies. In this paper, we place such static reductions into three categories: (i) those that are policy-independent (as those introduced by Witsenhausen in \cite{wit88}), (ii) those that are policy-dependent (as those introduced by Ho and Chu \cite{HoChu, ho1973equivalence} for partially nested dynamic teams), and (iii) those that we will refer to as static measurements with control-sharing reduction (where the measurements are static although control actions are shared according to the partially nested information structure). For the first type, we show that there is a bijection between person-by-person optimal (globally optimal) policies of dynamic teams and their policy-independent static reductions. For the second type, although there is a bijection between globally optimal policies of dynamic teams with partially nested information structures and their static reductions, in general there is no bijection between person-by-person optimal policies of dynamic teams and their policy-dependent static reductions. We also establish a stronger negative result concerning stationary solutions. We present sufficient conditions under which bijection relationships hold. Under static measurements with control-sharing reduction, connections between optimality concepts can be established under relaxed conditions. An implication is a convexity characterization of dynamic team problems under static measurements with control-sharing reduction. Finally, we consider multi-stage teams where we introduce equivalent models under which a single {\it agent} acting over the horizon is a collection of decision makers with increasing information over time (unlike the intrinsic model of Witsenhausen) and by taking into account an agent-wise optimality concept, we introduce two classes of ``agent-wise" static reductions: (i) {\it independent data} reduction under which the policy-independent reduction holds through agents and time, and (ii) {\it agent-wise (partially) nested independent reduction} under which measurements are independent through agents but (partially) nested through time for each agent. We study similar problems as that of single-stage setup for multi-stage problems, and we show that although there is a bijection between agent-wise person-by-person optimal (globally optimal) policies under both classes of reductions, there is no bijection between one-shot decision maker-wise (as considered earlier) person-by-person optimal policies in general under the {nested} reduction. Several illustrative examples are studied in detail.  Part II of the paper addresses similar issues in the context of stochastic dynamic games, where further subtleties arise.
\end{abstract}


%

\section{Introduction}

\label{sec:intro}

Team problems entail a collection of decision makers (DMs) acting together to optimize a common cost function, but not necessarily sharing all the available information. 
At each time stage, each DM has only partial access to the global information, which is characterized by the \textit{information structure} of the problem  \cite{wit75}. If there is a pre-defined order in which the DMs act, then the team is called a \textit{sequential team}. For sequential teams, if each DM's information depends only on primitive random variables, the team is \textit{static}. If at least one DM's information is affected by an action of another DM, the team is said to be \textit{dynamic}. Information structures can be further categorized as \textit{classical}, \textit{partially nested} (or quasi-classical), and \textit{nonclassical}. An information structure is \textit{classical} if the information of decision maker $i$ (DM$^i$) includes all of the information available to DM$^k$ for $k < i$. An information structure is \textit{partially nested}, if whenever the action of DM$^k$, for some $k < i$, affects the information of DM$^i$, then the information of DM$^i$ includes the information of DM$^k$. An information structure that is not partially nested is \textit{nonclassical}. 

For teams with finitely many DMs, Marschak \cite{mar55} has studied static teams and Radner \cite{Radner} has established connections between person-by-person (pbp) optimality, stationarity, and global optimality. Radner's results were generalized in \cite{KraMar82} by relaxing optimality conditions. The essence of these results is that in the context of static team problems, convexity of the cost function, subject to minor regularity conditions, suffices for the global optimality of pbp optimal solutions. In the particular case of LQG (Linear Quadratic Gaussian) static teams, this result leads to optimality of linear policies \cite{Radner}.

Optimality of linear policies also holds for dynamic LQG teams with  partially nested information structures through a transformation of the dynamic team to a static one \cite{HoChu}: in \cite{HoChu} for dynamic LQG teams with partially nested information structures and in \cite{ho1973equivalence} for general dynamic teams with partially nested information structures, satisfying an invertibility assumption (see Assumption \ref{assump:inv} later in this paper), it has been shown that they can be reduced to static team problems, where the aforementioned results for static teams can be applied. The transformation of  dynamic teams to static teams is called {\it{static reduction}}. In the static reduction presented in  \cite{HoChu, ho1973equivalence}, given the policies of the DMs, there is a bijection between observations as a function of precedent actions of DMs and   the primitive random variables, and observations generated under the transformations, where now they are only functions of primitive random variables.  We note that the static reduction in  \cite{HoChu, ho1973equivalence} depends on the policies that precedent DMs choose, and hence, in this paper, these will be referred to as  {\it policy-dependent static reduction}s (see Section \ref{sec:policysta}). 

On the other hand, in \cite{wit88}, Witsenhausen has introduced a static reduction for dynamic teams, where observations satisfy an absolute continuity condition and the information structure can be nonclassical, classical or partially nested. In this static reduction, the probabilistic nature of the problem has been transformed to the cost function by changing the measures of the observations to fixed probability measures. Witsenhausen's static reduction is independent of the policies that precedent DMs choose, and hence, we refer to this type of static reduction as a {\it{policy-independent static reduction}} (see Section \ref{sec:2.2}). The policy-independent static reduction is essentially a version of Girsanov's transformation \cite{girsanov1960transforming, benevs1971existence} which has been considered first in \cite[Eqn(4.2)]{wit88}, and later utilized in \cite[p. 114]{YukselBasarBook} and \cite[Section 2.2]{YukselWitsenStandardArXiv} (for discrete-time partially observed stochastic control, similar arguments had been presented, e.g. by Borkar in \cite{Bor00}, \cite{Bor07}). We refer the reader to \cite{charalambous2016decentralized} for relations with the classical continuous-time stochastic control, where the relation with Girsanov's classical measure transformation \cite{girsanov1960transforming, benevs1971existence} is recognized.

 Since Witsenhausen’s paper \cite{wit88}, the static reduction method has been shown to be very effective in arriving at existence, structural and approximation results. For existence results building on this approach, we refer the reader to \cite{gupta2014existence,YukselSaldiSICON17,SaldiArXiv2017,YukselWitsenStandardArXiv}, for a dynamic programming formulation to \cite{WitsenStandard} for countable spaces and \cite{YukselWitsenStandardArXiv} for general spaces, for rigorous approximations with finite models to \cite{saldiyuksellinder2017finiteTeam}, and for games to Part II \cite{SBSYgamesstaticreduction2021}. 

In this paper, we study the connections between pbp optimal (stationary, globally optimal) policies for both types of static reductions. 

We note that in the language of stochastic control, two interpretations of policies for a dynamic team and its static reduction can be stated as follows. For deterministic optimal control problems, open-loop policies do not explicitly depend on the history of the process, and these can be viewed as policies under a static information; whereas closed-loop policies can be viewed as policies under a dynamic information structure. Likewise, for single-DM classical stochastic control, {\it{(path-dependent) feedback}} policies (policies that are functions of a subset of the history of states or a noisy observation of states, which may depend on actions of the precedent DMs) are policies under a dynamic information structure, whereas {\it{noise feedforward}} policies (policies that are functions of only disturbances) are policies with a static  information structure \cite{benevs1971existence}. A subtlety of these connections for stationary (pbp optimal) policies can stem from the following observation: for dynamic teams, deviating a policy of a DM and fixing policies of others, requires a multi-directional deviation analysis since observations of frozen DMs depend on the deviating DM's policy. However, under a static reduction since observations of frozen DMs do not depend on the actions of the deviating DMs, only considering a single-directional deviation analysis of control actions is sufficient. Therefore, in general, establishing such connections between optimal/stationary policies is non-trivial and can fail to hold even under a partially nested information structure. In this paper, we present negative results and also sufficient conditions for positive results on the connections between optimality concepts of dynamic teams and their policy-independent and policy-dependent static reductions.

\begin{figure}[t!]
\begin{centering}
\tikzstyle{place}=[rectangle,draw=black!50,fill=blue!10,thick,
inner sep=2pt,minimum size=8mm]
\begin{tikzpicture}[scale=0.9]
\node at ( -2,2) [place, text width=3.3cm] (1){Pbp optimal policy for ($\mathcal{P}^{CS}$)};
\node at ( -2,0) [place, text width=3.3cm] (2){Pbp optimal policy for ($\mathcal{P}^{S}$)};
\node at ( 5.2,0) [place, text width=3.3cm] (3){Pbp optimal policy for ($\mathcal{P}^{D}$)};
\node at ( 5.2,2) [place, text width=3.3cm] (4){Pbp optimal policy for ($\mathcal{P}^{D, CS}$)};
\node at ( -3,1) [circle, inner sep=0pt, text width=3cm, red] (12){Corollary \ref{the:csstationary policies}};
\node at (6.9, 1) [circle, inner sep=1pt, text width=3cm]  (5){Theorem \ref{lem:2}};
\node at (-0.2, 1) [circle, inner sep=1pt, text width=3cm]  (5){Theorem \ref{lem:2}};

\node at (1.9, 1.5) [circle, inner sep=1pt, text width=2.5cm, blue]  (5){Theorem \ref{lem:2}};

\draw [blue, <->] (1)  to (4);
\draw [->] (2) to (1);
\draw [->] (3) to (4);
\draw [red, ->] (1) [out=180,in=180] to (2);
\end{tikzpicture}
\caption{Diagram of the connections between three optimality concepts for dynamic teams and their static measurements with control-sharing reductions.}\label{fig:3}
\end{centering}
\end{figure}
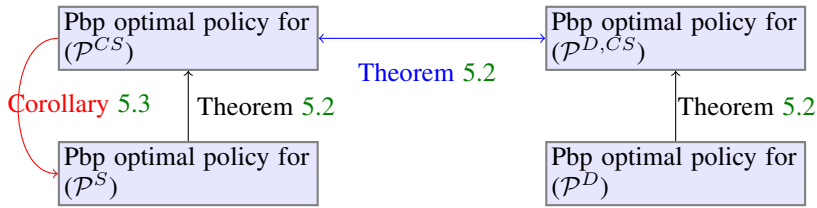

{\bf Contributions.} 
Our main contributions of this paper are summarized below:
\begin{itemize}[wide]
\item[(i)] For policy-dependent static reductions of stochastic teams with  partially nested information structure, in the reduced form, the cost functions are unaltered; this is not the case for the policy-independent static reduction of stochastic teams. On the other hand, the probability measure on the exogenous random variables do not change under the policy-dependent static reduction, but does so in the policy-independent static reduction. Furthermore, for the policy-dependent case, the static-reduced policy/map from exogenous variables to actions change depending on the policies that precedent DM choose, but in the policy-independent case, this map is unaltered (as the measurement variables are interpreted as exogenous variables). For policy-dependent static reductions, it is essential that policies are deterministic; however, policy-independent static reduction applies even when the policies are randomized.

\item[(ii)] We show in Theorem \ref{lem:1} that there is a bijection  between pbp optimal (globally optimal) policies (and under a further condition between stationary policies) of dynamic teams and their policy-independent static reductions (see Fig. \ref{fig:2.1}). This equivalency in relationships follows from the fact that this static reduction is policy-independent.

\item [(iii)] While for global optimality, policy-dependent static reductions and the dynamic information structure have equivalent optimal policies, for the policy-dependent case when one considers pbp optimal or stationary policies, significant subtleties emerge: a policy which is pbp optimal in one form may not be so in the other form (see Proposition \ref{the:negative}). Under sufficient convexity and minor regularity conditions on the cost function (see Assumption \ref{assump:2}) and a further regularity condition on policies and observations (see Condition (C)), we show in Theorem \ref{the:stationary policies} that there is a bijection between stationary (pbp optimal) policies of dynamic teams and their policy-dependent static reductions (see Fig. \ref{fig:2.2}). 

\item[(iv)] We define the reduction of dynamic stochastic teams with {\it partially nested with control-sharing} information structure to ones with {\it static measurements with control-sharing} information structure as static measurements with control-sharing reduction. We show that this reduction is independent of policies (see Theorems \ref{the:deptoind} and \ref{lem:2}), and facilitate our analysis in establishing the connections between optimality concepts (see Section \ref{sec:cs} for details and see Fig. \ref{fig:3} for a summary of the connections).

\item[(v)] The above static measurements with control-sharing reduction has implications on convexity properties: a partially nested dynamic team is convex in policies\footnote{see \cite[Section 3]{YukselSaldiSICON17} for the definition of convexity in policies.} if and only if its static measurements with control-sharing reduction is convex. This implies in general that for convexity (in policies), actions required to be shared under the policy-dependent static reduction requires (see Section \ref{sec:cvonx}).

\item [(vi)] Several examples are presented to illustrate both positive and negative results. Due to uniqueness of pbp optimal policies for LQG team problems under the policy-dependent static reduction, we establish stronger results for LQG models with a partially nested information structure (Corollary \ref{corollary:un1}).

\item [(vii)] We also study multi-stage team problems as a special setting of dynamic teams under two classes of static reductions: (i) {independent data} reduction under which the policy-independent reduction holds through agents and time, and (ii) {AG-wise (partially) nested independent reduction} under which measurements are independent through agents but (partially) nested through time. We show that there is a bijection between agent-wise pbp optimal policies (globally optimal policies) under both classes of reductions, but there is no bijection between (one-shot)-DM-wise pbp optimal policies in general under the {nested} reduction (Corollary \ref{corol:multi}). Furthermore, we discuss the impact of independent-data and AG-wise (partially) nested independent reductions on the variational analysis (Corollary \ref{multi-var}).  
\end{itemize}

{The organization of the paper is as follows: In Section \ref{sec:2.2}, we present preliminaries and provide a description of team problems within policy-independent, policy-dependent, and static measurements with control-sharing reductions. In Sections \ref{sec:varistatic}, \ref{sec:policydep}, and \ref{sec:cs}, we present results for dynamic teams under policy-independent, policy-dependent static reductions, and static measurements with control-sharing reduction, respectively. Multi-stage team problems are studied in Section \ref{sec:ms}. The paper ends with the concluding remarks of Section \ref{conc}, and several appendices.} This is Part I of a two-part paper, where Part II \cite{SBSYgamesstaticreduction2021} deals with stochastic dynamic games, presenting counterparts of the results in this paper for such multi-agent decision problems.

 \begin{figure}[t!]
\begin{centering}
\tikzstyle{place}=[rectangle,draw=black!50,fill=blue!10,thick,
inner sep=2pt,minimum size=8mm]
\begin{tikzpicture}[scale=0.9]
\node at ( 0,3) [place, text width=3.8cm] (1){Pbp optimal policy under static reduction};
\node at ( 5,3) [place,  text width=3.5cm] (2){Globally optimal policy under static reduction};
\node at ( -5,3) [place,  text width=3.5cm] (3){Stationary policy under static reduction};
\node at ( 0,0.9) [place, text width=3.5cm] (4){Pbp optimal policy for ($\mathcal{P}$)};
\node at ( 5,0.9) [place, text width=3.5cm](5){Globally optimal policy for ($\mathcal{P}$)};
\node at ( -5,0.9) [place, text width=3.5cm] (6){Stationary policy for ($\mathcal{P}$)};
\node at (2.8,2.5) [circle, inner sep=0pt, text width=1cm, red] (7){ \cite{YukselSaldiSICON17} };
\node at ( -2.3,2.7) [circle, inner sep=0pt, text width=1cm, red] (8){ \cite{YukselSaldiSICON17} };
\node at (-5.1,1.8) [circle, inner sep=0pt, text width=3.5cm, blue] (7){ Theorem \ref{lem:1} };
\node at (-0.2,1.8) [circle, inner sep=0pt, text width=3.5cm] (7){ Theorem \ref{lem:1} };
\node at (4.9,1.8) [circle, inner sep=0pt, text width=3.5cm] (7){ Theorem \ref{lem:1} };
\draw [<->] (1) to (4);
\draw [->] (2) to (1);
\draw [red, <->] (3) to (1);
\draw [<->] (2) to (5);
\draw [blue, <->] (3) to (6);
\draw [red, ->] (1) [out=-35,in=-155] to (2);
\draw [->] (5) to (4);
\draw [<->] (2) to (5);
\end{tikzpicture}
\caption{Diagram of the connections between three optimality concepts in dynamic teams and their policy-independent static reductions.}\label{fig:2.1}
\end{centering}
\end{figure}
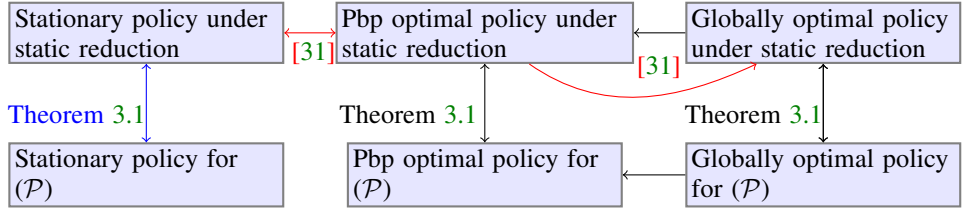

\begin{figure}[t!]
\begin{centering}
\tikzstyle{place}=[rectangle,draw=black!50,fill=blue!10,thick,
inner sep=2pt,minimum size=8mm]
\begin{tikzpicture}[scale=0.9]
\node at ( 0,3) [place, text width=3.5cm] (1){Pbp optimal policy for ($\mathcal{P}^{S}$)};
\node at ( 5,3) [place,  text width=3.5cm] (2){Globally optimal policy for ($\mathcal{P}^{S}$)};
\node at ( -5,3) [place,  text width=3.5cm] (9){Stationary policy for ($\mathcal{P}^{S}$)};
\node at ( 0,1) [place, text width=3.5cm] (4){Pbp optimal policy for ($\mathcal{P}^{D}$)};
\node at ( 5,1) [place, text width=3.5cm](3){Globally optimal policy for ($\mathcal{P}^{D}$)};
\node at ( -5,1) [place, text width=3.5cm] (20){Stationary policy for ($\mathcal{P}^{D}$)};
\node at (2.8,2.5) [circle, inner sep=0pt, text width=1cm, red] (70){ \cite{KraMar82}};
\node at ( -2.3,2.7) [circle, inner sep=0pt, text width=1cm, red] (80){ \cite{KraMar82}};
\node at (-5.7,2) [circle, inner sep=0pt, text width=3.5cm, blue] (90){Theorem \ref{the:stationary policies} };
\node at (-0.2,2) [circle, inner sep=0pt, text width=3.5cm, blue] (100){ Theorem \ref{the:stationary policies} };
\node at (4.8,2) [circle, inner sep=0pt, text width=3.5cm] (110){ Theorem \ref{the:4.2} };
\node at ( -4.1,2) [circle, red, inner sep=0pt, text width=2cm] (11){\Large$\times$};
\node at (-3,2) [circle, inner sep=0pt, text width=3cm] (60){Proposition \ref{the:negative}};

\draw [<->] (2) to (3);
\draw [blue, <->] (1)  to (4);
\draw [->] (2) to (1);
\draw [->] (3) to (4);
\draw [red, ->] (1) [out=-35,in=-155] to (2);
\draw [red, <->] (9) to (1);
\draw [<->] (9)  to (20);
\draw [blue, <->] (9) [out=180,in=180] to (20);
\end{tikzpicture}
\caption{Diagram of the connections between three optimality concepts for dynamic teams and their policy-dependent static reductions.}\label{fig:2.2}
\end{centering}
\end{figure}

\section{Information Structures and Static Reductions of Dynamic Teams}\label{sec:2.2}

\subsection{Witsenhausen's Intrinsic Model}\label{sec:pre}
Hans Witsenhausen's contributions \cite{wit75,WitsenhausenSIAM71,wit68,wit88,WitsenStandard} to stochastic control theory, and his characterization of information structures in decentralized stochastic control have been foundational in our modern understanding of decentralized stochastic control and decision theory. In this section, we introduce the characterizations as laid out by Witsenhausen, termed as {\it the Intrinsic Model} \cite{wit75}. In this model (described in discrete time), any action applied at any given time is regarded as applied by an individual DM, who acts only once. One particular advantage of this model, in addition to its generality, is that through it the characterizations regarding information structures can be concisely described. For a more comprehensive overview and some recent studies of information structures, we refer the reader to \cite{YukselBasarBook, CDCTutorial, saldiyukselGeoInfoStructure}.

Consider decentralized systems where DMs act in a pre-defined order.  Such systems are called {\it sequential teams}\footnote{see Andersland and Teneketzis \cite{AnderslandTeneketzisI}, \cite{AnderslandTeneketzisII} and Teneketzis \cite{Teneketzis2}, in addition to Witsenhausen \cite{WitsenhausenSIAM71} and \cite[p. 113]{YukselBasarBook} for non-sequential teams.}. For this class of teams, we now introduce Witsenhausen's \textit{intrinsic model} \cite{wit75}. 

\begin{itemize}[wide]
\item There exists a collection of \textit{measurable spaces} $\{(\Omega, {\cal F}), \allowbreak(\mathbb{U}^i,{\cal U}^i), (\mathbb{Y}^i,{\cal Y}^i), i \in {\mathcal{N}}\}$, specifying the system's distinguishable events, and control and measurement spaces. The set $\mathcal{N} :=\{1,2,\dots, N\}$ denotes the set of all DMs; the pair $(\Omega, {\cal F})$ is a
measurable space (on which an underlying probability may be defined); the pair $(\mathbb{U}^i, {\cal U}^i)$
denotes the measurable space from which the action $u^i$ of DM$^i$ is selected; and the pair $(\mathbb{Y}^i,{\cal Y}^i)$ denotes the measurable observation/measurement space. {Here, action and observation spaces for each DM are standard Borel spaces (that is, Borel subsets of complete, separable and metric spaces).}

\item There is a \textit{measurement constraint} to establish the connections between the observation variables and the system's distinguishable events. The $\mathbb{Y}^i$-valued observation variables are given by $y^i=h^i(\omega,{\underline u}^{[1,i-1]})$, where ${\underline u}^{[1,i-1]}=\{u^k~|~ k \leq i-1\}$ and $h^i$s are measurable functions.
\item The set of admissible control laws $\underline{\gamma}= \{\gamma^i\}_{i \in \mathcal{N}}$, also called
{\textit{designs}} or {\textit{policies}}, are measurable control functions, so that $u^i = \gamma^i(y^i)$. Let $\Gamma^i$ denote the set of all admissible policies for DM$^i$, and let ${\Gamma} = \prod_{i \in \mathcal{N}} \Gamma^i$.
\item There is a {\textit{probability measure}} $P$ on $(\Omega, {\cal F})$, making it a probability space on which the system is defined.
\end{itemize}

{For sequential teams, if each DM's information depends only on the primitive random variables, then the team is termed \textit{static}. If at least one DM's information is affected by an action of another DM, the team is said to be \textit{dynamic}. Information structures can be further categorized as \textit{classical}, \textit{partially nested}, and \textit{nonclassical}. An information structure is \textit{classical} if the information of DM$^i$ includes all of the information available to DM$^k$ for $k < i$. An information structure is \textit{partially nested} (or quasi-classical), if whenever the action of DM$^k$, for some $k < i$, affects the information of DM$^i$, then the information of DM$^i$ includes the information of DM$^k$. An information structure, which is not partially nested is \textit{nonclassical}. }

In view of Witsenhausen's equivalent model \cite{WitsenStandard}, any two information structures are equivalent if three conditions hold: (i) expected costs are the same, (ii) the measurable admissible policies are isomorphic (that is, measurability conditions enforced by the information structures are satisfied under equivalent models), and (iii) the constraints in the admissible policies are isomorphic. However, we will show that optimality properties of policies under equivalent information structures is rather fragile depending on the optimality concept considered, where an isomorphism of an optimal policy under an information structure may not be optimal under an equivalent information structure in general. In the next three subsections, we provide a description of dynamic teams with policy-independent and policy-dependent static reductions as well as static measurements with control-sharing reduction. Finally, in the last subsection, we develop refinements for a class of multi-stage problems. 

\subsection{Stochastic Dynamic Teams under Policy-Independent Static Reductions}\label{sec:staticred}

Let action and observation spaces be subsets of appropriate Euclidean spaces, i.e.,  $\mathbb{U}^{i} \subseteq \mathbb{R}^{n_{i}}$ and $\mathbb{Y}^{i} \subseteq \mathbb{R}^{m_{i}}$,  for $i \in \mathcal{N}$, where $n_{i}$ and $m_{i}$ are  positive integers. We formally introduce a stochastic dynamic team problem as follows:

\begin{itemize}[wide]
\item[\bf\text{Problem} \bf{($\mathcal{P}$)}:]
Consider a stochastic team problem within the intrinsic model with observations for each DM given by 
\begin{flalign}\label{eq:observy}
y^{i} = h^{i} (\omega_{0}, \omega_{i}, u^{1}, \dots, u^{i-1}, y^{1}, \dots, y^{i-1}),
\end{flalign}
where $\omega_{i}$ is an exogenous random variable, for $i \in \mathcal{N}$. Here,  $\omega_{0}$ is an $\Omega_{0}$-valued cost function-relevant exogenous random variable, i.e.,  $\omega_{0}:(\Omega,\mathcal{F}, P) \to (\Omega_{0},\mathcal{F}_{0})$, where $\Omega_{0}$ is a Borel space with its Borel $\sigma$-field $\mathcal{F}_{0}$. Let the information structure of DM$^{i}$ be $I^{i}=\{{y}^{i}\}$ (or $I^{i}=\{{y}^{k}\}_{k \in K_{i}}$ for a subset $K_{i}\subseteq \{1,\dots, i\}$). An expected cost function (to be minimized) under a policy $\underline{\gamma}= (\gamma^1, \cdots, \gamma^N) \in \Gamma$ is given by
\begin{equation}\label{eq:1.1}
\scalemath{0.94}{J(\underline{\gamma}) = E^{\underline{\gamma}}[c(\omega_{0},\underline{u})]:= E[c(\omega_{0},\gamma^1(y^1),\cdots,\gamma^N(y^N))]},
\end{equation}
for some Borel measurable cost function $c: \Omega_{0} \times \prod_{i=1}^{N} \mathbb{U}^i \to \mathbb{R}_{+}$ and $\underline {u}:=\{u^1, \cdots, u^{N}\}$.
\end{itemize}
\qedwhite

We first recall Witsenhausen's static reduction (see \cite{wit88, YukselWitsenStandardArXiv}), and then provide a description of team problems under this static reduction. Toward this end, we introduce an absolute continuity condition under which a policy-independent static reduction exists (see  \cite{YukselWitsenStandardArXiv} for further discussions).
\begin{assumption}\label{assump:pist}
For every $i\in \cal{N}$, there exists a probability measure $Q^{i}$ on $\mathbb{Y}^{i}$ and a function $f^{i}$ such that for any Borel set $A^{i}$
\begin{flalign}
&{P(y^{i} \in A^{i} \big|\omega_{0}, u^{1}, \dots, u^{i-1}, y^{1},\dots, y^{i-1})}\scalemath{0.94}{=\int_{A^i} f^{i}(y^{i},\omega_{0}, u^{1}, \dots, u^{i-1},y^{1},\dots, y^{i-1})Q^{i}(dy^{i})}\label{eq:abscon}.
\end{flalign}
\end{assumption}
Denote the joint distribution on $(\omega_{0}, u^{1},\dots, u^{N}, y^{1}, \dots, y^{N})$ by ${{P}}$, and the distribution of $\omega_{0}$ by $\mathbb{P}^{0}$. If the preceding absolute continuity condition holds, then there exists a joint reference distribution ${\mathbb{Q}}$ on $(\omega_{0}, u^{1},\dots, u^{N}, y^{1}, \dots, y^{N})$ such that the distribution ${{P}}$ is absolutely continuous with respect to ${\mathbb{Q}}$ (${{P}} \ll {\mathbb{Q}}$), where for every Borel set $A$ in $(\Omega_{0}\times \prod_{i=1}^{N}(\mathbb{U}^{i}\times\mathbb{Y}^{i}))$
\begin{flalign}
&\scalemath{0.94}{{{P}}(A)}\scalemath{0.94}{= \int_{A}\frac{d{{P}}}{d{\mathbb{Q}}}{\mathbb{Q}}(d\omega_{0}, du^{1},\dots, du^{N}, dy^{1}, \dots, dy^{N})}\label{eq:nonrandom},
\end{flalign}
where the joint reference distribution and Radon-Nikodym derivative (defined $P$-almost surely) are as follows:
\begin{flalign}
&{\mathbb{Q}}(d\omega_{0}, du^{1},\dots, du^{N}, dy^{1}, \dots, dy^{N}):=\mathbb{P}^{0}(d\omega_0)\prod_{i=1}^{N} Q^{i}(dy^{i}) 1_{\{\gamma^{i}(y^{i}) \in du^{i}\}},\label{eq:Q}\\
&\frac{d{{P}}}{d{\mathbb{Q}}}:=\frac{d{{P}}}{d{\mathbb{Q}}}(\omega_{0},u^{1}, \dots, u^{1},y^{1}\dots, y^{N}):=\prod_{i=1}^{N}f^{i}(y^{i},\omega_{0}, u^{1}, \dots, u^{i-1},y^{1},\dots, y^{i-1})\label{eq:dpdq}.
\end{flalign}

In view of the above derivations, we now formally introduce policy-independent static reductions:

\begin{definition}[{\bf Policy-Independent Static Reduction}]
For a stochastic team ($\mathcal{P}$) with cost function $c$ and a given information structure under Assumption \ref{assump:pist}, a policy-independent static reduction is defined as a change of measure \eqref{eq:nonrandom} under which measurements $y^{i}$ for each DM (defined in \eqref{eq:observy}) have independent distributions $Q^{i}$ and the expected cost function is given by
\begin{flalign}
{J(\underline{\gamma})}
&{:=E_{{\mathbb{Q}}}^{\underline{\gamma}}[\tilde{c}(\omega_0, u^{1},\dots, u^{N}, y^{1}, \dots, y^{N})]}\label{eq:2.4},
\end{flalign}
where the new cost function under the reduction is
\begin{flalign}
&{\tilde{c}(\omega_0, u^{1},\dots, u^{N}, y^{1}, \dots, y^{N})}{:= c(\omega_0, u^{1},\dots, u^{N}) \frac{d{{P}}}{d{\mathbb{Q}}}}.\label{eq:costpi10}
\end{flalign} \qedwhite
\end{definition}

\begin{remark}
\hfill 
\begin{itemize}[wide]
\item [(i)] As Witsenhausen noted, a static reduction always holds when the measurements take values from countable sets since a reference measure always exists on the measurement space $\mathbb{Y}^i$ (e.g., $Q^i(z) = \sum_{j \geq 1} 2^{-j} 1_{\{z = m_j\}}$, where $\mathbb{Y}^i=\{m_j~|~ j \in \mathbb{N}\}$) so that \eqref{eq:abscon} holds; 
\item [(ii)] The policy-independent static reduction applies even if policies are randomized. In fact, in \eqref{eq:Q}, we can replace the indicator function with a stochastic kernel $\Pi^{i}(u^{i}\in \cdot|y^{i})$ representing an independently randomized policy of DM$^{i}$ for $i \in \mathcal{N}$ (see also \cite[Section 2.2]{YukselSaldiSICON17}).
\item [(iii)] We emphasize that this static reduction is policy-independent since the change of measure \eqref{eq:abscon} is independent of policies of the precedent DMs. 
\end{itemize}
\end{remark}

We now recall definitions of globally optimal, pbp optimal, and stationary policies for ($\mathcal{P}$).
\begin{definition}[{\bf Optimality concepts for a dynamic team ($\mathcal{P}$)}]\label{eq:gof}
For a stochastic team ($\mathcal{P}$) with a given information structure, and cost function $c$:
\begin{itemize}[wide]
\item a policy \/ ${\underline \gamma}^{*}:=({\gamma^1}^*,\ldots, {\gamma^N}^*)\in { \Gamma}$\@ is
 \textit{globally optimal} if
\begin{equation*}
J({\underline \gamma}^*)=\inf_{{{\underline \gamma}}\in {{\Gamma}}}
J({{\underline \gamma}}):=\inf_{{{\underline \gamma}}\in {{\Gamma}}}E_{{{P}}}^{\underline{\gamma}}[{c}(\omega_0, u^{1},\dots, u^{N})] ,
\end{equation*}
\item a policy \/ ${\underline \gamma}^{*}\in { \Gamma}$\@ is \textit{pbp optimal (also called a Nash equilibrium)} if for all\/ $\beta
\in \Gamma^i$\@ and all\/ $i\in {\cal N}$\@,
\begin{equation*}\label{eq:1.3}
J({\underline \gamma}^*) \leq J({\underline \gamma}^{-i*},
\beta):=E_{{{P}}}^{({\underline \gamma}^{-i*},
\beta)}[{c}(\omega_0, u^{1},\dots, u^{N})],
\end{equation*}
where
$({\underline \gamma}^{-i*},\beta):= (\gamma^{1*},\ldots, \gamma^{(i-1)*},
\beta, \gamma^{(i+1)*},\ldots, \gamma^{N*})$,
\item a policy $\underline{\gamma}^{*}\in { \Gamma}$ is stationary if, for all $i \in \mathcal{N}$, $P$-a.s.,
\begin{flalign*}
&\nabla_{u^{i}}E_{{{P}}}\bigg[{c\bigg(\omega_{0},(\underline\gamma^{-i*}(y^{-i}), u^{i})\bigg) \bigg| y^{i}\bigg)\bigg]\bigg|_{u^{i}=\gamma^{i*}(y^{i})}=0},
\end{flalign*}
where $(\underline\gamma^{-i*}(y^{-i}), u^{i}):=(\gamma^{1*}(y^{1}),\dots, \gamma^{(i-1)*}(y^{i-1}), u^{i}, \gamma^{(i+1)*}(y^{i+1}),\dots, \gamma^{N*}(y^{N}))$.
\end{itemize}
\end{definition}
 \qedwhite

\begin{definition}[{\bf Optimality concepts under policy-independent static reduction}]\label{def:staticrde}
For a stochastic team ($\mathcal{P}$) with a given information structure, and cost function $\tilde{c}$ under a policy-independent static reduction:
\begin{itemize}[wide]
\item a policy \/ ${\underline \gamma}^{*}\in { \Gamma}$\@ is
 \textit{globally optimal} if
\begin{flalign*}
J({\underline \gamma}^*)=\inf_{{{\underline \gamma}}\in {{\Gamma}}}
J({{\underline \gamma}}):=\inf_{{{\underline \gamma}}\in {{\Gamma}}}E_{{\mathbb{Q}}}^{\underline{\gamma}}[\tilde{c}(\omega_0, u^{1},\dots, u^{N}, y^{1},\dots, y^{N})] ,
\end{flalign*}
\item a policy \/ ${\underline \gamma}^{*}\in { \Gamma}$\@ is \textit{pbp optimal}  if for all\/ $\beta^{i}
\in \Gamma^i$\@ and all\/ $i\in {\cal N}$\@,
\begin{flalign*}
J({\underline \gamma}^*) \leq J({\underline \gamma}^{-i*},
\beta):=E_{{\mathbb{Q}}}^{({\underline \gamma}^{-i*},
\beta)}[\tilde{c}(\omega_0, u^{1},\dots, u^{N}, y^{1},\dots, y^{N})],
\end{flalign*}
 \item a policy $\underline{\gamma}^{*}\in { \Gamma}$ is a stationary policy if $P$-a.s.,
\begin{flalign*}
&{\nabla_{u^{i}} E_{{\mathbb{Q}}}\bigg[\tilde{c}\bigg(\omega_{0},(\underline\gamma^{-i*}(y^{-i}), u^{i}),y^{1}, \dots, y^{N}\bigg) \bigg| y^{i}\bigg]\bigg|_{u^{i}=\gamma^{i*}(y^{i})}=0}.\nonumber
\end{flalign*}
\qedwhite
\end{itemize}
\end{definition}

One of our goals here is to study the connections between Definitions \ref{eq:gof} and \ref{def:staticrde}. In Section \ref{sec:varistatic}, we show the existence of a bijection between pbp optimal (globally optimal) policies of dynamic teams and their policy-independent static reductions, and  between stationary policies of dynamic teams and their policy-independent static reductions under a further condition on the Randon-Nikodym derivative (see \eqref{eq:extrastationary}). These connections are depicted in Fig \ref{fig:2.1}.

\subsection{Partially Nested Dynamic Teams under Policy-Dependent Static Reduction}\label{sec:policysta}
In the following, we first briefly recall Ho and Chu's static reduction \cite{HoChu, ho1973equivalence} (the policy-dependent static reduction), and then we provide a description of dynamic team problems under the policy-dependent static reduction. Consider stochastic dynamic teams with a partially nested information structure, where observations of DMs are given by
\begin{flalign}\label{eq:infody}
y_{i}^{D} := \bigg\{y^{D}_{\downarrow i}, \hat{y}^{D}_{i} :=g_{i}(h_{i}(\zeta), u^{D}_{\downarrow i})\bigg\},
\end{flalign}
where $\zeta:=\{\omega_{0},\dots, \omega_{N}\}$ denotes the set of all relevant random variables (corresponding to the uncertainty of the team associated with  the cost function and observations), and $g_{i}$ and $h_{i}$ are measurable functions. In the above, $y_{i}^{D}$ denotes a prespecified subset of collections of observations of DMs, specifying observations that are used in the construction of $u_{i}^{D}$. Also, $y^{D}_{\downarrow i}$ is a subset of a collection of observations of precedent DMs, all DM$^{j}$s for $j\in \mathcal{N}$, such that $\hat{y}^{D}_{i}$ is affected by the actions of DM$^{j}$ and $\{\downarrow i\}:=\{j\:|\: \hat{y}_{i}^{D}~~\text{is affected by}~~ u^{j}\}$. Let $I^{i}_{D}=\{{y}^{D}_{i}\}$ and the space of admissible policies under this information structure be given by
\begin{flalign}
{\Gamma^{D} :=\bigg\{\underline{\gamma}^{D}:=({\gamma}^{D}_{1},\dots,\gamma^{D}_{N})\: \bigg| \:u^{D}_{i} = \gamma^{D}_{i}(y_{i}^{D})\:\:\:{\forall}~~ i \in \mathcal{N}\bigg\}}\label{eq:gammaD}.
\end{flalign}
Under the above formulation, we introduce a class of dynamic team problems as follows: 
\begin{itemize}[wide]
\item[\bf\text{Problem} \bf{($\mathcal{P}^{D}$)}:]
For a stochastic team with information structure $I^{i}_{D}$ (with measurements $y_{i}^{D}$ defined in \eqref{eq:infody}) for all $i \in \mathcal{N}$, consider an expected cost function as in \eqref{eq:1.1} under the policy $\underline{\gamma}^{D}$. Derive a policy ${\underline \gamma}^{D*}:=({\gamma_1^{D*}},\dots, \gamma_N^{D*})\in \Gamma^{D}$\@ that is \textit{globally optimal} for ($\mathcal{P}^{D}$), that is
\begin{equation*}
J({\underline \gamma}^{D*})=\inf_{{{\underline \gamma}^{D}}\in {\Gamma}^{D}}
J({\underline \gamma}^{D}).
\end{equation*}  
Furthermore, derive a policy ${\underline \gamma}^{D*}$\@ that is \textit{pbp optimal} for ($\mathcal{P}^{D}$), that is
\begin{equation}
J({\underline \gamma}^{D*})=\inf_{{{\gamma}^{D}_{i}}\in {\Gamma}^{D}_{i}}
J({\gamma}^{D}_{i}, {\underline \gamma}^{D*}_{-i})\:\:\:\:\text{for all}\:\:\: i \in \mathcal{N}\label{eq:pbpteamd},
\end{equation}  
where ${\underline \gamma}^{D*}_{-i}= (\gamma^{D*}_{1},\dots, \gamma^{D*}_{i-1}, \gamma^{D*}_{i+1}, \dots, \gamma^{D*}_{N})$.\qedwhite
\end{itemize}

Now, we introduce an assumption under which the policy-dependent static reduction exists \cite{HoChu, ho1973equivalence}.
\begin{assumption}\label{assump:inv}
For all $i \in \mathcal {N}$ and for every fixed $u_{\downarrow i}^{D}$, the function {$g_{i}(\cdot, u^{D}_{\downarrow i}) : h_{i}(\zeta) \mapsto \hat{y}^{D}_{i}$} is invertible for all realizations of $\zeta$.
\end{assumption}

Following \cite{HoChu, ho1973equivalence}, under Assumption \ref{assump:inv}, given a policy $\underline\gamma^{D}$, the observations within the policy-dependent static reduction can be defined as follows:
\begin{flalign}
y_{i}^{S} = \bigg\{y^{S}_{\downarrow i}, \hat{y}^{S}_{i} :=h_{i}(\zeta)\bigg\},\label{eq:infost}
\end{flalign}
where $h_{i}(\zeta)=g_{i}^{-1}(\hat{y}^{D}_{i}, \gamma^{D}_{\downarrow i}(y_{\downarrow i}^{D}))$. Let the information structure of DM$^{i}$ be $I^{i}_{S}=\{{y}^{S}_{i}\}$ and
\begin{flalign}
\Gamma^{S} :=\bigg\{\underline{\gamma}^{S}=({\gamma}^{S}_{1},\dots,\gamma^{S}_{N})\: \bigg{|}\: u^{S}_{i} = \gamma^{S}_{i}(y_{i}^{S})\:\:\:{\forall}~~ i \in \mathcal{N}\bigg\}\label{eq:gammaS}.
\end{flalign}

A notable example is the LQG setting, studied by Ho and Chu, where via a static reduction in the sense above, optimality of linear policies for partially nested LQG teams has been established. For various examples of policy-dependent static reductions, we refer the reader to \cite{HoChu, ho1973equivalence}. We define team problems under the policy-dependent static reduction as follows:

\begin{itemize}[wide]
\item[\bf\text{Problem} \bf{($\mathcal{P}^{S}$)}:]
For a stochastic team with information structure $I^{i}
_{S}$ (with measurements $y_{i}^{S}$ defined in \eqref{eq:infost}) for all $i \in \mathcal{N}$, consider an expected cost function as in \eqref{eq:1.1} under policy $\underline{\gamma}^{S}$. Derive a policy ${\underline \gamma}^{S*}:=({\gamma_1^{S*}},\dots, \gamma_N^{S*})\in \Gamma^{S}$\@ that is \textit{globally optimal} for ($\mathcal{P}^{S}$): 
\begin{equation*}
J({\underline \gamma}^{S*})=\inf_{{{\underline \gamma}^{S}}\in {\Gamma}^{S}}
J({\underline \gamma}^{S}).
\end{equation*}  
Furthermore, derive a policy ${\underline \gamma}^{S*}$\@ that is \textit{pbp optimal} for ($\mathcal{P}^{S}$), that is
\begin{equation}
J({\underline \gamma}^{S*})=\inf_{{{\gamma}^{S}_{i}}\in {\Gamma}^{S}_{i}}
J({\gamma}^{S}_{i}, {\underline \gamma}^{S*}_{-i})\:\:\:\:\text{for all}\:\:\: i \in \mathcal{N}\label{eq:pbpteams}.
\end{equation} 
\end{itemize}
\qedwhite

\begin{definition}[{\bf Policy-Dependent Static Reduction}]
Consider a stochastic dynamic team ($\mathcal{P}^{D}$) with a given partially nested information structure $I^{i}
_{D}$, where Assumption \ref{assump:inv} holds. A policy-dependent static reduction is defined as the reduction of a stochastic dynamic team ($\mathcal{P}^{D}$) to a static one ($\mathcal{P}^{S}$) (which has an equivalent information structure $I^{i}_{S}$), where under the reduction, the cost function is unaltered and measurements are static, and for a given admissible policy ${\underline \gamma}^{D}\in \Gamma^{D}$, an admissible policy ${\underline \gamma}^{S}\in \Gamma^{S}$ can be constructed through a relation   
\begin{flalign}\label{eq:orpds}
u^{i} = \gamma^{S}_{i}(y_{i}^{S}) = \gamma^{D}_{i}(y_{i}^{D})~~\text{$P$-a.s.}
\end{flalign}
for all $i\in \cal{N}$. \qedwhite
\end{definition}

The relation \eqref{eq:orpds} in the construction of policies under the policy-dependent static reduction can be viewed as a composition of policies with a bijection $F_{\underline\gamma}^{i}:y^{D}_{i}\mapsto y^{S}_{i}$ with the inverse $(F_{\underline\gamma}^{i})^{-1}:y^{S}_{i}\mapsto y^{D}_{i}$ for all $i\in\cal{N}$, where the existence of this bijection follows from Assumption \ref{assump:inv} (we note that the subscript $\underline\gamma$ in the bijection denotes the fact that the bijection depends on the precedent policies for each $i\in \cal{N}$). That is, for any given policy $\underline\gamma^{D} \in \Gamma^{D}$, an admissible policy ${\underline \gamma}^{S}\in \Gamma^{S}$ can be constructed as $\gamma^{S}_{i}:=\gamma^{D}_{i} \circ F_{\gamma^{S}_{\downarrow i}}^{i}$, and for any given policy $\underline\gamma^{S}$, an admissible policy ${\underline \gamma}^{D}$ can be constructed as $\gamma^{D}_{i}:=\gamma^{S}_{i} \circ (F_{\gamma^{D}_{\downarrow i}}^{i})^{-1}$.

Some of our results in this paper address the following question: 
\begin{itemize}[wide]
\item [{\it Question 1:}] Given a stationary (pbp optimal, globally optimal) policy ${\underline \gamma}^{S*}\in \Gamma^{S}$ for ($\mathcal{P}^{S}$), is a policy ${\underline \gamma}^{D*}\in \Gamma^{D}$, constructed through relation \eqref{eq:orpds}, stationary (pbp optimal, globally optimal) policy for ($\mathcal{P}^{D}$)? Is the converse statement also true?
\end{itemize}
\begin{remark}
\hfill
\begin{itemize}[wide]
\item[(i)] In contrast to the policy-independent static reduction, the policy-dependent static reduction requires DMs to have access to (able to compute) the actions of precedent DMs according to a partially nested information structure; hence, it requires the policies to be deterministic. However, for teams, without any loss of optimality, globally optimal policies can be chosen among those that are deterministic \cite[Theorems 2.3 and 2.5]{YukselSaldiSICON17}.
\item[(ii)] We also note that in the policy-dependent static reduction, in contrast to the policy-independent static reduction, the cost function will not change under the static reduction.
\end{itemize}
\end{remark}

In Section \ref{sec:policydep}, we first show that the answer to Question 1 is affirmative for globally optimal policies of dynamic teams and their policy-dependent static reduction. However, for pbp optimal and stationary policies, there might not exist a bijection between stationary (pbp optimal) policies of dynamic teams and their policy-dependent static reduction in general (see Fig. \ref{fig:2.2}). This is quite opposite of the result for the policy-independent static reduction (see Theorem \ref{lem:1}), where there is a bijection between stationary (pbp optimal, globally optimal) policies of dynamic teams and their policy-independent static reduction. We also present sufficient conditions such that the answer to Question 1 is affirmative for pbp optimal and stationary policies. Several examples including LQG models are presented.

\subsection{Partially Nested with Control-Sharing Information Structure and Static Measurements with Control-Sharing Reduction}\label{sec:pncs}
{To establish connections between pbp optimality and convexity of dynamic teams and their policy-dependent static reductions, we introduce {\it (dynamic) partially nested with control-sharing} team problems, where we expand the information structure such that in addition to observations, actions are also shared (this expansion is consistent with partially nested information structure), i.e., for each DM$^{i}$,  
\begin{flalign}
y_{i}^{D, CS} := \bigg\{y^{D}_{\downarrow i}, u^{\downarrow i},\hat{y}^{D}_{i}\bigg\}\label{eq:cs-Dver},
\end{flalign}
with $I_{i}^{D, CS}:=\{y_{i}^{D, CS}\}$ and
\begin{flalign}
\Gamma^{D, CS} :=\bigg\{\underline{\gamma}^{D, CS}=({\gamma}^{D, CS}_{1},\dots,\gamma^{D, CS}_{N})\: \bigg{|}\: u^{D, CS}_{i} = \gamma^{D, CS}_{i}(y_{i}^{D, CS})\:\:\:{\forall}~~ i \in \mathcal{N}\bigg\}\label{eq:gammaDCS}.
\end{flalign}

\begin{itemize}[wide]
\item[\bf\text{Problem} \bf{($\mathcal{P}^{D, CS}$)}:]
For a stochastic team with information structure $I^{D, CS}_{i}$ (with measurements \eqref{eq:cs-Dver}), consider the expected cost function (to be minimized) as in \eqref{eq:1.1} under policy $\underline{\gamma}^{D, CS}$. 
\qedwhite
\end{itemize}

Under the invertibility condition (Assumption \ref{assump:inv}), there is a bijection between dynamic observations $y_{i}^{D}$ and static one $y_{i}^{S}$ for each $i \in \mathcal{N}$, and hence, this allows us to reduce the original dynamic team to another one where measurements are static. We refer to this reduction as {\it static measurements with control-sharing}. We define the observations within this reduction as follows 
\begin{flalign}
&y^{CS}_{i}:=\bigg\{y^{S}_{\downarrow i}, u^{\downarrow i}, \hat{y}^{S}_{i} \bigg\}\label{eq:control-sharing information structure},
\end{flalign}
with $I^{CS}_{i}:=\{y^{CS}_{i}\}$ and
\begin{flalign}
\Gamma^{CS} :=\bigg\{\underline{\gamma}^{CS}=({\gamma}^{CS}_{1},\dots,\gamma^{CS}_{N})\: \bigg{|}\: u^{CS}_{i} = \gamma^{CS}_{i}(y_{i}^{CS})\:\:\:{\forall}~~ i \in \mathcal{N}\bigg\}\label{eq:gammaCS}.
\end{flalign}

\begin{itemize}[wide]
\item[\bf\text{Problem} \bf{($\mathcal{P}^{CS}$)}:]
For a stochastic team with information structure $I^{CS}_{i}$ (with measurements \eqref{eq:control-sharing information structure}), consider an expected cost function (to be minimized) as in \eqref{eq:1.1} under policy $\underline{\gamma}^{CS}$.\qedwhite
\end{itemize}

We refer to the above problems as {\it static measurements with control-sharing} team problems: 
\begin{definition} [{\bf Static Measurements with Control-Sharing Reduction}] Consider a stochastic dynamic partially nested with control-sharing stochastic team ($\mathcal{P}^{D, CS}$) with a given information structure $I^{D, CS}_{i}$, where Assumption \ref{assump:inv} holds. A static measurements with control-sharing reduction is defined as the reduction of a stochastic dynamic team ($\mathcal{P}^{D, CS}$) to a static measurements control-sharing problem ($\mathcal{P}^{CS}$) with information structure, $I^{CS}_{i}$, where under the reduction, the cost function is unaltered and the measurements are static (see \eqref{eq:control-sharing information structure}), and for a given admissible policy $\underline{\gamma}^{D, CS}$ for ($\mathcal{P}^{D, CS}$), an admissible policy $\underline{\gamma}^{CS}$ for ($\mathcal{P}^{CS}$) can be constructed, through the relation  
\begin{flalign}
{\gamma}^{D, CS}_{i}(y_{i}^{D, CS}) = \gamma^{CS}_{i}(y_{i}^{CS})~~\text{for every $u^{\downarrow i}$ $P$-a.s.}\label{eq:stmr}
\end{flalign} 
for all $i\in \cal{N}$. \qedwhite
\end{definition}

In Section \ref{sec:policydep}, based on static measurements with control-sharing reduction, we establish some isomorphic connections between pbp optimal (globally optimal, stationary) policies ($\mathcal{P}^{D}$), ($\mathcal{P}^{S}$), ($\mathcal{P}^{CS}$), and ($\mathcal{P}^{D, CS}$). We also study convexity of team problems under the policy-dependent and static measurements with control-sharing reductions (see Fig. \ref{fig:3}).}

\subsection{Multi-Stage Team Problems and Static Reductions}\label{sec:multi}

We now consider multi-stage stochastic dynamic teams and introduce two further reductions in the contexts of policy-independent and policy-dependent static reductions introduced for single-stage team problems in the preceding subsections. We first recall that under the intrinsic model of Witsenhausen (see Section \ref{sec:intro}), every DM acts separately and once (which we refer to as the {\it one-shot-DM property} in the following discussion). However, depending on the information structure and the cost function, it may be convenient to consider a collection of DMs as {\it a single agent} acting at different time instants. In fact, in classical stochastic control, this is the standard approach. With this motivation, we will introduce a new reduction concept building on the one introduced by Witsenhausen (called independent-data reduction) \cite[Section 2.4]{wit88} and another one in \cite[Section 3.2]{sanjari2019optimal}. The underlying idea is to view DMs acting in a sequence with increasing information as a single agent with a larger action space. This facilitates our optimality analysis. We note that this approach, in a less general form, was utilized to establish structural and existence results in  \cite[Section 3.2]{sanjari2019optimal}. 

\begin{itemize}[wide]
\item[\bf\text{Problem} \bf{($\mathcal{P}^{\text{Multi}}$)}:]
Consider the following formulation of multi-stage stochastic teams:
\item[(i)] The state dynamics and observations for $t \in \mathcal{T}:=\{0,\dots, T-1\}$ are given respectively by 
\begin{flalign}
x_{t+1} &= f_{t}(x_{0:t}, u_{0:t}^{1:N}, w_{t}),\\
y_{t}^{i} &= h_{t}^{i}(x_{0:t}, u_{0:t-1}^{1:N}, v_{t}^{i})\label{eq:msments},
\end{flalign}
for all $i\in \mathcal{N}$, where $f_{t}$ and $h_{t}^{i}$ are measurable functions. $x_{0:t}:=(x_{0}, \dots, x_{t})$, and $w_{t}, v_{t}^{1}, \dots, v_{t}^{N}$ for all $t \in\mathcal{T}$ are  random variables taking values in standard Borel spaces. Further, we let $u_{0:t}^{1:N}:=(u_{0}^{1}, \dots, u_{t-1}^{1}, \dots, u_{0}^{N}, \dots, u_{t}^{N})$, and introduce appropriate collections of DMs as agents, with the $i$-th agent (AG$^{i}$) for $i\in \mathcal{N}$ acting at different time instants $t\in \mathcal{T}$ and comprised of DM$^{i}_{0}$, $\dots$, DM$^{i}_{T-1}$. 
\item[(ii)]  The observation, action, state, and disturbance spaces are standard Borel spaces with ${\bf{Y}}^{i}:=\prod_{t=0}^{T-1}\mathbb{Y}^{i}_{t}$, ${\bf{U}}^{i}:=\prod_{t=0}^{T-1}\mathbb{U}^{i}_{t}$, ${\bf{X}}:=\prod_{t=0}^{T-1}\mathbb{X}_{t}$, ${\bf{W}}^{i}:=\prod_{t=0}^{T-1}\mathbb{W}^{i}_{t}$, ${\bf{V}}^{i}:=\prod_{t=0}^{T-1}\mathbb{V}^{i}_{t}$, respectively. 
\item[(iii)] An admissible policy for AG$^{i}$ is $\pmb{\gamma^{i}}\in {\bf {\Gamma}^{i}}$, where $\pmb{\gamma^{i}}:=(\gamma^{i}_{0}, \dots,\gamma^{i}_{T-1})$ and  ${\bf {\Gamma^{i}}}=\prod_{t=0}^{T-1}\Gamma^{i}_{t}$. Given an information structure $I_{t}^{i} \subseteq \{y_{0:t}^{1:N}, u_{0:t-1}^{1:N}\}$, each admissible policy $\gamma_{t}^{i}$ is a measurable function with $u_{t}^{i}=\gamma^{i}_{t}(I_{t}^{i})$. 
\item [(iv)] A multi-stage expected cost function under a policy $\underline{\pmb{\gamma}}$  is given by
\begin{flalign}\label{eq:multi1.1}
&J(\underline{\pmb{\gamma}})={E}^{\underline{\pmb{\gamma}}}\bigg[\sum_{t=0}^{T-1}c_{t}(\omega_{0}, x_{t}, u_{t}^{1},\dots, u_{t}^{N})+c_{T}(x_{T})\bigg],
\end{flalign}
for some Borel measurable cost function $c: \Omega_{0} \times \mathbb{X}_{t} \times \prod_{i=1}^{N} \mathbb{U}_{t}^{i} \to \mathbb{R}_{+}$, where $\underline{\pmb{\gamma}}:=(\pmb{\gamma^{1}},\pmb{\gamma^{2}},\dots,\pmb{\gamma^{N}})$, and again $\omega_{0}$ is an $\Omega_{0}$-valued cost function-relevant exogenous random variable, $\omega_{0}:(\Omega,\mathcal{F}, {P}) \to (\Omega_{0},\mathcal{F}_{0})$, where $\Omega_{0}$ is a Borel space with its Borel $\sigma$-field $\mathcal{F}_{0}$.
\end{itemize}

Determine (existence and characterization of) the policy $\underline{\pmb{\gamma}}$ that minimizes \eqref{eq:multi1.1} over $\prod_{i=1}^{N}{\bf {\Gamma^{i}}}$.\qedwhite

\begin{definition}\label{def:plwise}
For a multi-stage stochastic team, a policy $\underline{\pmb{\gamma}}^{*}$ is agent-wise (AG-wise) pbp optimal if for all $i\in \mathcal{N}$ and for all ${\pmb{\beta}} \in {\bf{\Gamma^{i}}}$,
\begin{flalign*}
J(\underline{\pmb{\gamma}}^*) \leq J({\pmb{\gamma}^{-i,*}}, {\pmb{\beta}}),
\end{flalign*}
where ${\pmb{\gamma}^{-i*}}:=(\pmb{\gamma}^{1*}, \dots, \pmb{\gamma}^{i-1*},\pmb{\gamma}^{i+1*},\dots,\pmb{\gamma}^{N*})$. A policy $\underline{\pmb{\gamma}}^{*}$ is AG-wise globally optimal if for all $\underline{\pmb{\gamma}} \in \prod_{i=1}^{N}{\bf{\Gamma^{i}}}$,
\begin{flalign*}
J(\underline{\pmb{\gamma}}^*) \leq J(\underline{\pmb{\gamma}}). 
\end{flalign*}
Also, a policy $\underline{\pmb{\gamma}}^{*}$ is (one-shot) DM-wise pbp optimal if for all $i\in \mathcal{N}$ and $k \in \mathcal{T}$ and for all ${{\beta}^{i}_{t}} \in {{\Gamma^{i}_{t}}}$,
\begin{flalign*}
J(\underline{\pmb{\gamma}}^*) \leq J({\pmb{\gamma}^{-i*}}, (\gamma^{i*}_{-t}, {{\beta}^{i}_{t}})),
\end{flalign*}
where $(\gamma^{i*}_{-t}, {{\beta}^{i}_{t}}):=(\gamma^{i*}_{0}, \dots, \gamma^{i*}_{t-1}, {{\beta}^{i}_{t}}, \gamma^{i*}_{t+1}, \dots, \gamma^{i*}_{T-1})$. \qedwhite

\end{definition}
 
Comparing Definitions \ref{def:plwise} and \ref{eq:gof}, we can see that concepts of AG-wise and DM-wise global optimality are equivalent. Also, every AG-wise pbp optimal policy is DM-wise pbp optimal; however, the converse statement is not true in general. The reason is that in the definition of AG-wise pbp optimality, in contrast to the definition of DM-wise pbp optimality, policies $(\gamma_{0}^{i}, \dots, \gamma_{t-1}^{i}, \gamma_{t+1}^{i}, \dots, \gamma_{T-1}^{i})$ are not frozen. Later on, we will provide sufficient conditions for the converse statement to hold based on our static reductions for multi-stage teams. However, we note that clearly if there is a unique DM-wise pbp optimal policy, then it is a unique AG-wise pbp optimal policy. In the following, we first introduce two assumptions used for our static reductions, and introduce static reductions for multi-stage stochastic teams.

\begin{assumption}\label{assump:indatapolicy}
For every $i\in \mathcal{N}$ and every $t\in \mathcal{T}$, there exists a probability measure $\tilde{Q}^{i}_{t}$ on $\mathbb{Y}_{t}^{i}$ and a measurable function $\phi^{i}_{t}$ such that for all Borel sets $A=A^{1}\times \dots \times A^{N}$ with $A^{i}$ in $\mathbb{Y}_{t}^{i}$, we have 
\begin{flalign}
&{P\bigg((y_{t}^{1},\dots,y_{t}^{N})} \scalemath{0.94}{\in A ~\bigg| ~\omega_{0}, x_{0}, v_{0:t-1}^{1:N}, w_{0:t-1}^{1:N}, y^{1:N}_{0:t-1}, u_{0:t-1}^{1:N}\bigg)}\nonumber\\
&{\:\:\:\:\:\:\:=\prod_{i=1}^{N}\int_{A^i}\phi_{t}^{i}(y^{i}_{t}, \omega_{0},x_{0}, v_{0:t-1}^{1:N}, w_{0:t-1}^{1:N},y_{0:t-1}^{1:N}, u_{0:t-1}^{1:N})\tilde{Q}_{t}^{i}(dy_{t}^{i})}.\label{eq:indepdata}
\end{flalign}
\end{assumption}

Let $\tilde{\mathbb{P}}$ be the joint distribution on $(\omega_{0}, {x}_{0}, \underline{\pmb{w}},\underline{\pmb{v}},\underline{\pmb{u}}, \underline{\pmb{y}})$, and $\mu$ be the fixed joint distribution on $(\omega_{0}, {x}_{0}, \underline{\pmb{w}},\underline{\pmb{v}})$. Let $\underline{\pmb{z}}:=(\pmb{z}^{1},\dots, \pmb{z}^{N})$ and $\pmb{z}^{i}:=(z_{0}^{i},\dots, z_{T-1}^{i})$ for $z=u, y, w, v$ and $i\in \mathcal{N}$. Hence, under the preceding change of measure \eqref{eq:indepdata}, there exists a joint reference distribution $\tilde{\mathbb{Q}}$ on $(\omega_{0}, {x}_{0}, \underline{\pmb{w}},\underline{\pmb{v}},\underline{\pmb{u}}, \underline{\pmb{y}})$ such that $\tilde{\mathbb{P}}$ is absolutely continuous with respect to $\tilde{\mathbb{Q}}$, where for every Borel set $B$ on $(\Omega_{0}\times \mathbb{X}_{0}\times \prod_{i=1}^{N}({\bf{W}}^{i} \times {\bf{V}}^{i} \times {\bf{U}}^{i} \times {\bf{Y}}^{i}))$
\begin{flalign}
&\tilde{\mathbb{P}}(B)= \int_{B}\frac{d{\tilde{\mathbb{P}}}}{d\tilde{\mathbb{Q}}}{\tilde{\mathbb{Q}}}(d\omega_{0}, d{x}_{0}, d\underline{\pmb{w}},d\underline{\pmb{v}},d\underline{\pmb{u}}, d\underline{\pmb{y}})\label{eq:nonrandom3},
\end{flalign}
where
\begin{flalign*}
\tilde{\mathbb{Q}}(d\omega_{0}, d{x}_{0}, d\underline{\pmb{w}},d\underline{\pmb{v}},d\underline{\pmb{u}}, d\underline{\pmb{y}})&:=\mu(d\omega_{0}, d{x}_{0}, d\underline{\pmb{w}},d\underline{\pmb{v}})\prod_{t=0}^{T-1}\prod_{i=1}^{N}\tilde{Q}_{t}^{i}(dy_{t}^{i})1_{\{\gamma^{i}_{t}(y^{i}_{t}) \in du^{i}_{t}\}},\\
\frac{d{\tilde{\mathbb{P}}}}{d\tilde{\mathbb{Q}}}&:=\prod_{t=0}^{T-1}\prod_{i=1}^{N}\phi_{t}^{i}(y^{i}_{t}, \omega_{0},x_{0}, v_{0:t-1}^{1:N}, w_{0:t-1}^{1:N},y_{0:t-1}^{1:N}, u_{0:t-1}^{1:N}),
\end{flalign*}
where $\mu$ is the fixed distribution on $(\omega_{0}, {x}_{0}, \underline{\pmb{w}},\underline{\pmb{v}})$.
\begin{assumption}\label{assum:decop}
For every $i \in \mathcal{N}$, there exists a probability measure $\hat{Q}^{i}$ such that for every Borel set $B$
\begin{flalign}
\tilde{\mathbb{P}}(B) &= \int_{B} \frac{d\tilde{\mathbb{P}}}{d\hat{\mathbb{Q}}} \hat{\mathbb{Q}}(d\underline{\pmb{u}},d\underline{\pmb{y}}, d\underline{\pmb{w}}, d\omega_0)\label{eq:chng} ,\\
\hat{\mathbb{Q}}(d\underline{\pmb{u}},d\underline{\pmb{y}}, d\underline{\pmb{w}}, d\omega_0)&:=\prod_{i=1}^{N}\hat{Q}^{i}(d\pmb{u}^i,d\pmb{y}^i,d{\pmb{w}}^{i}) \mathbb{P}^{0}(d\omega_{0})\nonumber.
\end{flalign}
We note that in the above, distributions $\tilde{\mathbb{P}}$, $\tilde{\mathbb{Q}}$, and  $\hat{Q}^{i}$ depend on policies; however, a change of measure \eqref{eq:chng} is policy-independent. 
\end{assumption}

\begin{definition}[{\bf Independent-Data and AG-wise (Partially) Nested Independent Reductions}]
Consider a multi-stage stochastic team ($\mathcal{P}^{\text{Multi}}$) with a given information structure. Introduce the following two agent-wise reductions for it:
\begin{itemize}[wide]
\item [(i)] (Independent-data reduction) 
Let Assumption \ref{assump:indatapolicy} hold. An independent-data reduction is a change of measure \eqref{eq:nonrandom3} under which the measurements driven by \eqref{eq:msments} have distributions $\tilde{Q}_{t}^{i}$, and the expected cost function can be written as follows:
\begin{flalign}
&J(\underline{\pmb{\gamma}}):=E^{\underline{\pmb{\gamma}}}_{\tilde{\mathbb{P}}}\bigg[\sum_{t=0}^{T-1}c_{t}(\omega_0, x_{t}, {u}^{1}_{t},\dots, {u}^{N}_{t})+c_{T}(x_{T})\bigg]=E^{\underline{\pmb{\gamma}}}_{\tilde{\mathbb{Q}}} \bigg[\hat{c}(\omega_0, {x}_{0}, \underline{\pmb{w}},\underline{\pmb{v}},\underline{\pmb{u}}, \underline{\pmb{y}})\bigg]\label{eq:2.4nested},
\end{flalign}
where the new cost function is
\begin{flalign}
&\hat{c}(\omega_0, {x}_{0}, \underline{\pmb{w}},\underline{\pmb{v}},\underline{\pmb{u}}, \underline{\pmb{y}}):= \sum_{t=0}^{T-1}c_{t}(\omega_0, x_{t}, {u}^{1}_{t},\dots, {u}^{N}_{t})\frac{d{\tilde{\mathbb{P}}}}{d\tilde{\mathbb{Q}}}.\label{eq:ncostms}
\end{flalign}
The team problem under this static reduction can be viewed as the one that Witsenhausen referred to as \textit{a static problem with independent data} \cite{wit88};
\item [(ii)] (AG-wise (partially) nested independent reduction) Let Assumption \ref{assum:decop} hold. AG-wise nested independent reduction is a reduction under which for each AG$^{i}$ through $t\in \mathcal{T}$, the information structure is nested (i.e., $\sigma(y_{t}^{i}) \subset \sigma(y_{t+1}^{i})$), and the expected cost function can be written as follows:
\begin{flalign*}
&J(\underline{\pmb{\gamma}})=E^{\underline{\pmb{\gamma}}}_{\hat{\mathbb{Q}}} \bigg[c(\omega_{0}, \underline{\pmb{u}},\underline{\pmb{y}}, \underline{\pmb{w}})\frac{d{\tilde{\mathbb{P}}}}{d\hat{\mathbb{Q}}}\bigg].
\end{flalign*}
If for each AG$^{i}$ through $t\in \mathcal{T}$, the information structure is only partially nested, the reduction is called an {\it AG-wise partially nested independent reduction}.

\end{itemize}
\qedwhite
\end{definition}


We note that one scenario where the AG-wise (partially) nested independent reduction arises is when each agent has a nested private information structure and the policy-independent reduction can be applied through agents (or only through dynamics and not necessarily for observation through time) such that under the reduction, Assumption \ref{assum:decop} holds. We also note that the independent-data reduction does not require the information structure to be nested, and on the other hand, the AG-wise (partially) nested independent reduction does not require Assumption \ref{assump:indatapolicy} to hold (see Examples \ref{lem:exmple} and \ref{lem:exmple1}). In particular, the AG-wise (partially) nested independent reduction can be applied even in the presence of common noise (or common random shocks to all agents through states or dynamics) without any further assumptions on the noise processes or the structures of the dynamics and observations. Furthermore, the AG-wise (partially) nested independent reduction also allows noiseless control and/or state sharing through time for each agent (where $y_{t}^{i} = h_{t}^{i}(x_{0:t}^{i}, u_{0:t-1}^{i})$). Later on, in Section \ref{sec:ms} (see Corollary \ref{corol:multi}), we show that AG-wise pbp optimal policies for (multi-stage) dynamic teams remain AG-wise pbp optimal policies for the teams under independent-data and AG-wise (partially) nested independent reductions; however, DM-wise pbp optimal policies only remain DM-wise pbp optimal policies for teams under independent-data static reductions and not necessarily under AG-wise (partially) nested independent reductions.

\begin{remark}
Two settings, where AG-wise based reductions are useful, are as follows:
\begin{itemize}[wide]
\item [(i)] Mean-field teams can be viewed as limit models of symmetric finite agent teams with a mean-field interaction (for example, see \cite{sanjari2019optimal, SSYdefinetti2020} for mean-field teams, and \cite{carmona2018probabilistic} and references therein for mean-field games). We note that for multi-stage mean-field dynamic teams, the independent-data and AG-wise nested independent reductions have been introduced in \cite[Section 3.2]{sanjari2019optimal} and \cite[Assumption 5.1(ii)]{SSYdefinetti2020}. As it has been shown in \cite[Section 3.2]{sanjari2019optimal} and \cite[Assumption 5.1(ii)]{SSYdefinetti2020}, the above static reduction  under mild conditions on the action and observation spaces leads to closedness of a set of policies for each agent through times under an appropriate topology, which is desirable for establishing existence and/or convergence results. 
\item [(ii)] We also note that the infinite horizon team problem under a AG-wise (partially) nested independent reduction is more tractable compared to an independent-data static reduction. That is because  AG-wise (partially) nested independent reductions allow agents to have nested information structures without requiring independent-data (which can be viewed as the total recall property of the private history for agents, where measurements may not necessarily be independent random variables under the reduction). Furthermore, using AG-wise (partially) nested independent reductions leads to richness in the variational analysis since for multi-stage team problems joint perturbations through times of a given agent are allowed (see Corollary \ref{multi-var} and Remark \ref{rem:NR}).
\end{itemize}
\end{remark}

\section{Optimal Policies for Dynamic Teams under Policy-Independent Static Reduction}\label{sec:varistatic}

In this section, we establish connections between stationary (pbp optimal, globally optimal) policies for dynamic teams and their policy-independent static reductions.

\begin{theorem}\label{lem:1}
Consider a stochastic dynamic team ($\mathcal{P}$) with a policy-independent static reduction \eqref{eq:abscon}.
\begin{itemize}[wide]
\item [(i)] A policy $\underline{\gamma}^*$ is pbp optimal (globally optimal)   for ($\mathcal{P}$) if and only if $\underline{\gamma}^*$ is pbp optimal (globally optimal) for a policy-independent static reduction of ($\mathcal{P}$);
\item [(ii)] Let a policy $\underline{\gamma}^*$ satisfy $P$-a.s.,
\begin{flalign}
&\nabla_{u^{i}}E^{\gamma^{-i*}}_{{\mathbb{Q}}}\bigg[\frac{d{{P}}}{d{\mathbb{Q}}}\bigg|y^{i}\bigg]\bigg|_{u^{i}=\gamma^{i*}(y^{i})}=0\:\:\:\:\forall i\in\cal{N}\label{eq:extrastationary},
\end{flalign} 
where $\frac{d{{P}}}{d{\mathbb{Q}}}$ is defined in \eqref{eq:dpdq}. Then, $\underline{\gamma}^*$ is stationary for ($\mathcal{P}$)  if and only if $\underline{\gamma}^*$ is stationary for a policy-independent static reduction of ($\mathcal{P}$).
\end{itemize}
\end{theorem}
\begin{proof}
Proof is provided in the Appendix.
\end{proof}

We note that \eqref{eq:extrastationary} implies that at a stationary point $\underline{\gamma}^*$, the decision laws have no impact locally on the reduction (on the Radon-Nikodym derivative term $\frac{d{{P}}}{d{\mathbb{Q}}}$). Now, in view of Theorem \ref{lem:1}, we use \cite[Theorem 3.3 and 3.4]{YukselSaldiSICON17} (which is a generalization of  \cite[Theorem 2 and 3]{KraMar82} by using an information structure dependent nature of convexity under policy-independent static reductions) to introduce sufficient conditions for stationary policies of ($\mathcal{P}$) to be globally optimal for dynamic teams ($\mathcal{P}$) with a policy-independent static reduction. We note that  if the team is static, then \eqref{eq:extrastationary} trivially holds since $\frac{d{{P}}}{d{\mathbb{Q}}}$ does not depend on actions, and hence, the following corollary to Theorem \ref{lem:1} provides a refinement for \cite[Theorems 3.3 and 3.4]{YukselSaldiSICON17} and \cite[Theorems 2 and 3]{KraMar82}.

\begin{corollary}\label{claim:1}
Consider a stochastic dynamic team ($\mathcal{P}$) with a policy-independent static reduction \eqref{eq:abscon}.  Assume that
\begin{itemize}[wide]
\item [(i)] The cost function, $c$, and the Radon-Nikodym derivative, $\frac{d{{P}}}{d{\mathbb{Q}}}$, are continuously differentiable in $u^{1}, \dots, u^{N}$,
\item [(ii)] The cost function under the policy-independent static reduction, $\tilde{c}$ (defined in \eqref{eq:costpi10}), is convex in $u^{1}, \dots, u^{N}$. 
\end{itemize}
Suppose that $\underline{\gamma}^{*}$ is a stationary policy for ($\mathcal{P}$) and satisfies \eqref{eq:extrastationary}. Let for all $\underline{\gamma} \in \Gamma$ with $E^{\underline{\gamma}}_{{{P}}} [{c}(\cdot)]<\infty$, the following two conditions hold for all $i\in \cal{N}$: 
\begin{flalign}
&E_{{{P}}}\bigg[\nabla_{u^{i}} {c}(\omega_{0},\underline{\gamma}^{*}(\underline{y}))\bigg(\gamma^{i}(y^{i}) -\gamma^{i*}(y^{i})\bigg)\bigg] <\infty,\label{eq:sum}\\
&E_{{\mathbb{Q}}}\bigg[\bigg(\nabla_{u^{i}}\frac{d{{P}}}{d{\mathbb{Q}}}(\omega_{0}, \underline{y},\underline{\gamma}^{*}(\underline{y}))\bigg)c(\omega_{0}, \underline{\gamma}^{*}(\underline{y}))\bigg(\gamma^{i}(y^{i}) -\gamma^{i*}(y^{i})\bigg)\bigg]<\infty,\label{eq:sumed1}
\end{flalign}
where 
\begin{flalign*}
&\nabla_{u^{i}}{c}(\omega_{0},\underline{\gamma}^{*}(\underline{y})):=\nabla_{u^{i}} {c}\bigg(\omega_{0},(\underline\gamma^{-i*}(y^{-i}), u^{i})\bigg)\bigg|_{u^{i}=\gamma^{i*}(y^{i})},\\
&\nabla_{u^{i}}\frac{d{{P}}}{d{\mathbb{Q}}}(\omega_{0}, \underline{y},\underline{\gamma}^{*}(\underline{y})):=\nabla_{u^{i}}\frac{d{{P}}}{d{\mathbb{Q}}}\bigg(\omega_{0}, \underline{y},(\underline\gamma^{-i*}(y^{-i}), u^{i})\bigg) \bigg|_{u^{i}=\gamma^{i*}(y^{i})}.
\end{flalign*}
Then, $\underline{\gamma}^{*}$ is globally optimal for ($\mathcal{P}$). Moreover, if the cost function $\tilde{c}$ is strictly convex in $(u^{1},\dots, u^{N})$, $\underline{\gamma}^{*}$ is the unique globally optimal policy for ($\mathcal{P}$).
\end{corollary}
\begin{proof}
The proof follows from Theorem \ref{lem:1}(ii)  and \cite[Theorem 3.3 and 3.4]{YukselSaldiSICON17}.
\end{proof}

\begin{remark}\label{remcvx}
\hfill
\begin{itemize}[wide]
\item [(i)] We note that for dynamic teams with a given information structure, if $c$ is convex in $u^{1}, \dots, u^{N}$ for all $\omega_{0}, \omega_{1}, \dots, \omega_{N}$, then $\tilde{c}$ (see \eqref{eq:costpi10}) is not necessarily convex in $u^{1}, \dots, u^{N}$ for all $\omega_{0}, \omega_{1}, \dots, \omega_{N}$. In particular, the celebrated Witsenhausen's counterexample \cite{wit68} is an example of non-convexity becoming evident under a policy-independent static reduction (this has been precisely shown in \cite[equation (3.5)]{ YukselSaldiSICON17}).
\item  [(ii)] We also note the result in \cite[Theorem 3.6]{YukselSaldiSICON17}, where it has been shown that a dynamic team with a policy-independent static reduction is convex in policies (see, \cite[Definition 3.1]{YukselSaldiSICON17}) if and only if its policy-independent static reduction is.
\item [(iii)] In view of Corollary \ref{claim:1}, even if the cost function $\tilde{c}$ in \eqref{eq:costpi10} is assumed to be convex and continuously differentiable in actions under policy-independent static reductions, \eqref{eq:sum}-\eqref{eq:sumed1} might not be sufficient to establish global optimality of a stationary policy of $\underline{\gamma}^{*}$ for ($\mathcal{P}$), in general. The reason is that, the effect of the deviating policies (or a deviating policy in the definition of the stationary policy $u^{i}$) on the probability measures of observations has not been taken into account. This observation suggests that in dynamic teams, variational analysis requires to take into account the effect of the deviating policies on the probability measures of observations (which has been considered in  \eqref{eq:extrastationary} and \eqref{eq:sumed1}). 
\end{itemize}
\end{remark}

\section{Optimal Policies for Dynamic Teams under Policy-Dependent Static Reductions}\label{sec:policydep}

Here, we study the connections between stationary (pbp optimal, globally optimal) policies of dynamic teams and their policy-dependent static reductions (see Section \ref{sec:policysta}). We first have the following result.
\begin{theorem}\label{the:4.2}
Consider a stochastic dynamic team ($\mathcal{P}^{D}$) with partially nested information structure. Let Assumption \ref{assump:inv} hold. Then, ${\underline \gamma}^{D*}$ is a globally optimal policy for ($\mathcal{P}^{D}$) if and only if ${\underline \gamma}^{S*}$ is a globally optimal policy for ($\mathcal{P}^{S}$) under the policy-dependent static reduction (see \eqref{eq:orpds}).
\end{theorem}
{\begin{proof}
Since the information structure is partially nested and policies are deterministic,  under Assumption \ref{assump:inv}, there is a bijection from the set of  policies $\Gamma^{D}$ to the set of  policies $\Gamma^{S}$. Therefore, global optimality in one domain implies global optimality in the other domain.  
\end{proof}}

\subsection{Stationary and PBP Optimal Policies for Dynamic Teams and their Policy-Dependent Static Reductions}\label{sec:stteam}

Here, we provide three examples that serve to demonstrate the subtlety of the connections between stationary (pbp optimal) policies of ($\mathcal{P}^{D}$) and ($\mathcal{P}^{S}$). These are counterexamples which show that, in contrast to the case of globally optimal policies, the isomorphism relations between stationary (pbp optimal) policies of ($\mathcal{P}^{D}$) and ($\mathcal{P}^{S}$) are no longer true, in general (under Assumption \ref{assump:inv}). 

We first show that a policy ${\underline \gamma}^{S, *}$ is  stationary (also pbp optimal) for ($\mathcal{P}^{S}$), but ${\underline \gamma}^{D*}$, satisfying the policy-dependent static reduction is not pbp optimal for ($\mathcal{P}^{D}$).

\begin{example}\label{ex:st2}
  Consider a $2$-DM stochastic team ($\mathcal{P}^{D}$) with $I^{1}=\{y^{D}_{1}\}$ and $I^{2}:=\{y^{D}_{2}\}=\{y^{D}_{1}, \hat{y}^{D}_{2}\}$, where $\hat{y}^{D}_{2}=\omega_{2}+{u^{1}}$, and $\omega_{2}=: \hat{y}^{S}_{2}$ and $y^{D}_{1}=: y^{S}_{1}=\omega_{1}$ are primitive random variables. Let the expected cost function be given as 
 \begin{flalign}
 E[c(\omega_{2},u^{1}, u^{2})] := E[(u^{1}-u^{2}+\omega_{2})^{2} - \alpha (u^{1})^{2}]\label{eq:costex1}, 
 \end{flalign}
 for a given $\alpha \in (0, 1)$.
 \begin{itemize}[wide]
 \item A policy $\underline{\gamma}^{S*}=({\gamma}^{S*}_{1}, {\gamma}^{S*}_{2})=(0, (0, I))$ (where the policy $(0, (0, I))$ denotes ${\gamma}^{S*}_{1} \equiv 0$, ${\gamma}^{S, *}_{2, 1} \equiv 0$, and ${\gamma}^{S*}_{2, 2}$ is the identity map, $I$, that is, $u^{1*}=\gamma^{S*}_{1}(y^{S}_{1})=0$ and $u^{2*}={\gamma}^{S*}_{2}(y^{S}_{1}, \hat{y}^{S}_{2})=\hat{y}^{S}_{2}$) is pbp optimal for ($\mathcal{P}^{S}$). 
 \item However, a policy $\underline{\gamma}^{D*}=({\gamma}^{D*}_{1}, {\gamma}^{D*}_{2})=(0, (-\gamma^{S*}_{1}, I))$ constructed under a relation \eqref{eq:orpds} (where the policy $(0, (-\gamma^{S*}_{1}, I))$ denotes ${\gamma}^{D*}_{1} \equiv 0$, ${\gamma}^{D*}_{2, 1}=-\gamma^{S*}_{1}$, and ${\gamma}^{D*}_{2, 2}$ is the identity map, that is, $u^{1*}=\gamma^{D*}_{1}(y^{D}_{1})=0$ and $u^{2*}=\hat{y}^{D}_{2}-\gamma^{S*}_{1}(y^{D}_{1})$) is not pbp optimal for ($\mathcal{P}^{D}$) since fixing a policy of DM$^{2}$ to $\gamma^{D*}_{2}$ such that $u^{2*}=\hat{y}^{D}_{2}-\gamma^{S*}_{1}(y^{D}_{1})$, the expected cost function will be concave in $u^{1}$ ($c(u^{1}, u^{2*})=-\alpha (u^{1})^{2}$) and the value will be unbounded from below. We note, however, that $\underline{\gamma}^{D*}$ is a stationary policy for ($\mathcal{P}^{D}$). 
 \end{itemize}
\end{example}

In the following example, we show that a policy ${\underline \gamma}^{D*}$ is  stationary (also pbp optimal) for ($\mathcal{P}^{D}$), but ${\underline \gamma}^{S, *}$ under the policy-dependent static reduction, is not pbp optimal for ($\mathcal{P}^{S}$).

\begin{example}\label{ex:st2inv}
  Consider a $2$-DM stochastic team ($\mathcal{P}^{D}$) with $I^{1}=\{y^{D}_{1}\}$ and $I^{2}:=\{y^{D}_{2}\}=\{y^{D}_{1}, \hat{y}^{D}_{2}\}$, where $\hat{y}^{D}_{2}=\omega_{2}+{u^{1}}$, and $\omega_{2}=: \hat{y}^{S}_{2}=\omega_{1}$ and $y^{D}_{1}=: y^{S}_{1}$ are primitive random variables. Let the expected cost function be given as 
 \begin{flalign}
 E[c(\omega_{2},u^{1}, u^{2})] := E[\alpha (u^{1})^{2}+ \beta (u^{2}-\omega_{2})^{2}-(u^{1}-u^{2}+\omega_{2})^{2}]\label{eq:costex1inv}, 
 \end{flalign}
 for a given $\alpha \in (0, 1)$ and $\beta >1$.
 \begin{itemize}[wide]
 \item A policy  $\underline{\gamma}^{D*}=({\gamma}^{D*}_{1}, {\gamma}^{D*}_{2})=(0, (0, I))$  (where ${\gamma}^{D*}_{1} \equiv 0$, ${\gamma}^{D*}_{2, 1} \equiv 0$ and ${\gamma}^{D*}_{2, 2}$ is the identity map, that is, $u^{1*}=\gamma^{D*}_{1}(y^{D}_{1})=0$ and $u^{2*}=\hat{y}^{D}_{2}$) is pbp optimal for ($\mathcal{P}^{D}$) since fixing the policy of DM$^{2}$ to ${\gamma}^{D*}_{2}$, the expected cost function will be convex in $u^{1}$ ($c(u^{1}, u^{2})=(\alpha+\beta) (u^{1})^{2}$), and fixing the policy of DM$^{1}$ to ${\gamma}^{D*}_{1}$ such that $u^{1}=\gamma^{D*}_{1}(y^{D}_{1})=0$, the expected cost function will be convex in $u^{2}$ ($c(u^{1}, u^{2})=(\beta-1) (u^{2}-\omega_{2})^{2}$). 
  \item However, under the policy-dependent static reduction, the policy $\underline{\gamma}^{S*}=({\gamma}^{S*}_{1}, {\gamma}^{S*}_{2})=(0, (-\gamma^{D*}_{1}, I))$ constructed under a relation \eqref{eq:orpds}, is not pbp optimal for ($\mathcal{P}^{S}$) since fixing the policy of DM$^{2}$ to ${\gamma}^{S*}_{2}$ such that $u^{2}={\gamma}^{S*}_{2}(y^{S}_{1}, \hat{y}^{S}_{2})=\hat{y}^{S}_{2}-\gamma^{D*}_{1}(y^{S}_{1})$, the expected cost function will be concave in $u^{1}$ ($c(u^{1}, u^{2})=(\alpha-1) (u^{1})^{2}$). 
 \end{itemize}
\end{example}

Next, we provide an example where a policy ${\underline \gamma}^{D, *}$ is  stationary (pbp optimal) for ($\mathcal{P}^{D}$), but the corresponding policy ${\underline \gamma}^{S, *}$ under the policy-dependent static reduction,  is not stationary for ($\mathcal{P}^{S}$).
\begin{example}\label{ex:st1}
 Consider a $2$-DM stochastic team ($\mathcal{P}^{D}$) with $I^{1}=\{y^{D}_{1}\}$ and $I^{2}=\{y^{D}_{2}\}:=\{y^{D}_{1}, \hat{y}^{D}_{2}\}$, where $\hat{y}^{D}_{2}=\omega_{2}+\sqrt{u^{1}}$, and $\omega_{2}$ and $y^{D}_{1}=y^{S}_{1}:=\omega_{1}$ are primitive random variables. Let $\mathbb{U}^{1}=\mathbb{R}_{+}$ and the expected cost function be given by 
 \begin{flalign}
 E[c(\omega_{2},u^{1}, u^{2})] := E[(\sqrt{u^{1}}-u^{2}+\omega_{2})^{2}].\label{eq:costex2}
 \end{flalign}
\begin{itemize}[wide]
\item A policy $\underline{\gamma}^{D*}=({\gamma}^{D*}_{1}, {\gamma}^{D*}_{2})=(0,(0, I))$ (where ${\gamma}^{D*}_{1} \equiv 0$, ${\gamma}^{D*}_{2, 1} \equiv 0$ and ${\gamma}^{D*}_{2, 2}$ is the identity map, that is, $u^{1*}=0$ and $u^{2*}=\hat{y}^{D}_{2}$) is stationary for ($\mathcal{P}^{D}$). 
\item However, under the policy-dependent static reduction, the corresponding policy $\underline{\gamma}^{S*}=({\gamma}^{S*}_{1}, {\gamma}^{S*}_{2})=(0, (\sqrt{\gamma^{D*}_{1}}, I))$ constructed under the relation \eqref{eq:orpds} (where ${\gamma}^{S*}_{1} \equiv 0$, ${\gamma}^{S*}_{2, 1}=\sqrt{\gamma^{D*}_{1}}$, and ${\gamma}^{S*}_{2, 2}$ is the identity map, that is, $u^{1}=0$ and $u^{2}=\omega_{2}+\sqrt{\gamma^{D*}_{1}(y^{S}_{1})}$) is not stationary (although it is pbp optimal) for ($\mathcal{P}^{S}$). Since fixing the policy of DM$^{2}$ to ${\gamma}^{S*}_{2}$ such that $u^{2}=\omega_{2}$, the derivative of the expected cost function with respect to $u^{1}$ is always $1$. Hence, the criterion for stationarity does not lead to a solution. 
\end{itemize}
\end{example}

Hence, in view of the preceding examples, we have the following negative result. 
\begin{prop}\label{the:negative}
Consider a stochastic dynamic team ($\mathcal{P}^{D}$) with partially nested information structure. Let Assumption \ref{assump:inv} hold. Then:
\begin{itemize}[wide]
\item [(i)] If ${\underline \gamma}^{D*}$ is stationary (pbp optimal) for ($\mathcal{P}^{D}$), then ${\underline \gamma}^{S*}$ is not necessarily  stationary (pbp optimal) for ($\mathcal{P}^{S}$) under the policy-dependent static reduction (see \eqref{eq:orpds});
\item[(ii)] If ${\underline \gamma}^{S*}$ is a stationary (pbp optimal) policy for ($\mathcal{P}^{S}$), then ${\underline \gamma}^{D*}$, satisfying the policy-dependent static reduction relation \eqref{eq:orpds}, is not necessarily pbp optimal for ($\mathcal{P}^{D}$).
\end{itemize}
\end{prop}
\begin{proof}
This is a direct consequence of the examples above, where Examples \ref{ex:st2inv} and \ref{ex:st1} imply Part (i), and Example \ref{ex:st2} implies Part (ii).
\end{proof}

Next, we introduce sufficient conditions to establish connections between stationary policies of dynamic teams and their policy-dependent static reductions.  We first introduce a regularity and convexity condition on the cost function and a regularity condition on observations and policies needed for our result (see Theorem \ref{the:stationary policies}).
\begin{assumption}\label{assump:2}
For every $\omega_{0}$,
\begin{itemize}
\item [(a)] the cost function $c$ is continuously differentiable in $(u^{1},\dots, u^{N})$, 
\item [(b)] the cost function $c$ is (jointly) convex in $(u^{1},\dots, u^{N})$.
\end{itemize}
\end{assumption}

\begin{condition1}
A policy ${\underline\gamma}^{D}$ satisfies Condition (C) if for all $i\in \mathcal{N}$, $P$-a.s, $\gamma^{D}_{i}(\{g_{j}(h_{j}(\zeta),u^{\downarrow j})\}_{j\in \downarrow i}, g_{i}(h_{i}(\zeta), u^{\downarrow i}))$ is affine in $u^{\downarrow i}$.
\end{condition1}

\begin{theorem}\label{the:stationary policies}
Consider a stochastic dynamic team ($\mathcal{P}^{D}$) with partially nested information structure. Let Assumptions \ref{assump:inv} and \ref{assump:2} hold. Then 
a policy ${\underline \gamma}^{D*}$ satisfying Condition (C) is a stationary (pbp optimal) policy for ($\mathcal{P}^{D}$) if and only if  ${\underline \gamma}^{S*}$ is a stationary (pbp optimal) policy for ($\mathcal{P}^{S}$) under the policy-dependent static reduction (see \eqref{eq:orpds}).  
\end{theorem}
\begin{proof}
Proof is provided in the Appendix.
\end{proof}


\section{Optimality and Convexity under Static Measurements with Control-Sharing Reduction}\label{sec:cs}
In this section, we present our results for the static measurements with control-sharing reduction and its impact on optimality and convexity of dynamic team problems.

\subsection{Optimal Policies under Static Measurements with Control-Sharing Reduction}
In this subsection, we consider team problems with partially nested with control-sharing information structure (see Section \ref{sec:pncs}), and  establish  connections between pbp optimal (globally optimal, stationary) policies of ($\mathcal{P}^{D}$), ($\mathcal{P}^{S}$), ($\mathcal{P}^{D, CS}$), and ($\mathcal{P}^{CS}$). 

Now, we state the following result as a theorem, since it will be consequential later on:
\begin{theorem}\label{the:deptoind}
For a stochastic dynamic team with a partially nested information structure, where Assumption \ref{assump:inv} holds, static measurements with control-sharing reduction is policy-independent.
\end{theorem}

\begin{proof}
Since Assumption \ref{assump:inv} holds and DMs have access to $u^{\downarrow i}$, static measurements with control-sharing reduction to ($\mathcal{P}^{CS}$) for each DM is independent of precedent DMs' policies: Given $\underline\gamma^{D, CS}$, a policy $\underline\gamma^{CS}$ can be constructed through \eqref{eq:stmr}, i.e., for every $i\in \mathcal{N}$, $u^{i}=\gamma^{D, CS}_{i}(y^{D}_{\downarrow i}, u^{\downarrow i},g_{i}(h_{i}(\zeta), u^{\downarrow i}))=\gamma^{CS}_{i}(y^{S}_{\downarrow i}, u^{\downarrow i},\hat{y}^{S}_{i})$ for every $u^{\downarrow i}$ $P$-a.s. The fact that the expected cost function does not change under the above reduction completes the proof.
\end{proof}

Since it is possible to construct examples in the spirit of Examples \ref{ex:st2}, \ref{ex:st2inv}, and \ref{ex:st1}, the isomorphic connections between stationary (pbp optimal) policies of ($\mathcal{P}^{CS}$) (($\mathcal{P}^{D, CS}$)) and ($\mathcal{P}^{S}$) and/or ($\mathcal{P}^{D}$) fail to hold in general. Now, we provide some positive results.

\begin{theorem}\label{lem:2}
Consider stochastic dynamic teams ($\mathcal{P}^{D}$), ($\mathcal{P}^{D, CS}$), and ($\mathcal{P}^{CS}$) with partially nested information structure.
\begin{itemize}[wide]
\item [(i)] If Assumption \ref{assump:inv} holds, then a policy $\underline{\gamma}^{D, CS}$ is pbp optimal (stationary, globally optimal) for ($\mathcal{P}^{D, CS}$) if and only if $\underline{\gamma}^{CS}$ is a pbp optimal (stationary, globally optimal) policy for ($\mathcal{P}^{CS}$) under the static measurements with control-sharing reduction (see \eqref{eq:stmr}).
\item [(ii)] Any pbp optimal (stationary) policy ${\underline \gamma}^{D*}\in \Gamma^{D}$  for ($\mathcal{P}^{D}$) constitutes a pbp optimal (stationary) policy on the enlarged space $\Gamma^{D, CS}$ for ($\mathcal{P}^{D, CS}$); however, in general, if ${\underline \gamma}^{D, CS*}\in \Gamma^{D, CS}$ is pbp optimal (stationary) for ($\mathcal{P}^{D, CS}$), then ${\underline \gamma}^{D, *}$ satisfying $\gamma^{D*}_{i}(y^{D}_{i})=\gamma^{D, CS*}_{i}(y^{D, CS}_{i})$ $P$-a.s. for all $i \in \mathcal{N}$, is not necessarily pbp optimal (stationary) for ($\mathcal{P}^{D}$).
\item [(iii)] Any pbp optimal (stationary) policy ${\underline \gamma}^{S*}\in \Gamma^{S}$ for ($\mathcal{P}^{S}$) constitutes a pbp optimal (stationary) policy on the enlarged space $\Gamma^{CS}$ for ($\mathcal{P}^{CS}$); however, in general, if ${\underline \gamma}^{CS*}\in \Gamma^{CS}$ is pbp optimal (stationary) for ($\mathcal{P}^{CS}$), then ${\underline \gamma}^{S*}$ satisfying $\gamma^{S*}_{i}(y^{S}_{i})=\gamma^{CS*}_{i}(y^{CS}_{i})$ $P$-a.s. for all $i \in \mathcal{N}$, is not necessarily pbp optimal (stationary) for ($\mathcal{P}^{S}$).
\end{itemize}
\end{theorem}
\begin{proof}
Proof is provided in the Appendix.
\end{proof}


Now, we state a corollary to Theorems \ref{the:stationary policies}, \ref{the:deptoind}, and \ref{lem:2}.

\begin{corollary}\label{the:csstationary policies}
Consider a stochastic dynamic team ($\mathcal{P}^{D}$) with partially nested information structure, where  Assumption \ref{assump:inv} holds.
\begin{itemize}[wide]
\item [(i)] A policy ${\underline \gamma}^{CS*}$ is globally optimal for ($\mathcal{P}^{CS}$) if and only if policies ${\underline \gamma}^{D*}$ and ${\underline \gamma}^{S*}$ are globally optimal for ($\mathcal{P}^{D}$) and ($\mathcal{P}^{S}$), respectively, with a relation $\gamma^{D*}_{i}(y^{D}_{i})=\gamma^{S*}_{i}(y^{S}_{i})=\gamma^{CS*}_{i}(y^{CS}_{i})$ $P$-a.s. for all $i \in \mathcal{N}$ (for any static representation of $\underline\gamma^{S*}$, there may exist multiple representations for $\underline\gamma^{D*}$ and $\underline\gamma^{CS*}$). 
\item [(ii)] Under  Assumption \ref{assump:2}(a), if a pbp optimal policy ${\underline \gamma}^{CS*}$ for ($\mathcal{P}^{CS}$) is affine in actions, then ${\underline \gamma}^{S*}$, satisfying $\gamma^{S*}_{i}(y^{S}_{i})=\gamma^{CS*}_{i}(y^{CS}_{i})$ $P$-a.s. for all $i \in \mathcal{N}$, is pbp optimal for ($\mathcal{P}^{S}$).
\end{itemize}
\end{corollary}
\begin{proof}
Part (i) follows from the fact that the expected cost function is identical for policies  $\underline{\gamma}^{D*}$, $\underline{\gamma}^{S*}$ and $\underline{\gamma}^{CS*}$, and Part (ii) follows from the fact that $\underline{\gamma}^{CS*}$ is a affine function of actions of precedent DMs, and hence, a similar argument as in the proof of Theorem \ref{the:stationary policies} completes the proof.
\end{proof}

\begin{remark}\label{rem:OP-CL}
Let $\tilde{\Gamma}^{CS}$ be a space of admissible policies \eqref{eq:gammaCS}, where observations \eqref{eq:control-sharing information structure} are replaced by 
\begin{flalign*}
\tilde{y}^{CS}_{i}:=\bigg\{y^{S}_{\downarrow i}, u^{K_{i}}, \hat{y}^{S}_{i} :=h^{i}(\zeta) \bigg\}
\end{flalign*}
for $K_{i} \subseteq {\downarrow i}$, and $\tilde{I}^{CS}=\{\tilde{y}^{CS}_{i}\}$, that is, for each $i \in \mathcal{N}$, control action of DM$^{i}$ is only shared to a subset $K_{i} \subseteq {\downarrow i}$ of DM$^{j}$s with ${j \in \downarrow i}$ (similarly, we can define $\tilde{\Gamma}^{D, CS}$ as \eqref{eq:gammaDCS} with expanded observations \eqref{eq:cs-Dver}). Following from the proof of Theorem \ref{lem:2}, the results of Theorem \ref{lem:2} (ii)(iii) remain valid if the enlarged space $\Gamma^{D, CS}$ is replaced by $\tilde{\Gamma}^{D, CS}$, where $\Gamma^{D} \subseteq \tilde{\Gamma}^{D, CS} \subseteq  \Gamma^{D, CS}$ (or the enlarged space $\Gamma^{CS}$ is replaced by $\tilde{\Gamma}^{CS}$ where $\Gamma^{S} \subseteq \tilde{\Gamma}^{CS} \subseteq  \Gamma^{CS}$) (we note that sets of admissible policies   $\Gamma^{D, CS},  \Gamma^{D},  \Gamma^{CS}$ and $\Gamma^{S}$ are introduced in  \eqref{eq:gammaDCS}, \eqref{eq:gammaD}, \eqref{eq:gammaCS}, and \eqref{eq:gammaS}, respectively).
\end{remark}

Next, we present results on the existence and uniqueness of optimal policies for ($\mathcal{P}^{D}$), ($\mathcal{P}^{S}$), ($\mathcal{P}^{D, CS}$), and ($\mathcal{P}^{CS}$), using Theorems \ref{the:4.2}, Corollary \ref{the:csstationary policies}, and Examples \ref{ex:st2} and \ref{ex:st2inv}.
\begin{prop}\label{corollary:un}
Consider a stochastic dynamic team ($\mathcal{P}^{D}$) with partially nested information structure, where Assumption \ref{assump:inv} holds. Then: 
\begin{itemize}[wide]
 \item [(i)] The existence of a pbp optimal policy for ($\mathcal{P}^{CS}$) (($\mathcal{P}^{D, CS}$)) does not imply the existence of a pbp optimal policy for ($\mathcal{P}^{S}$) (($\mathcal{P}^{D}$)).
 \item[(ii)]  The globally optimal policy is unique for ($\mathcal{P}^{S}$) if and only if the globally optimal policy is essentially unique for ($\mathcal{P}^{D}$); 
 \item [(iii)] If globally optimal policies for ($\mathcal{P}^{S}$) and/or ($\mathcal{P}^{D}$) are essentially unique, then there exists an essentially  unique globally optimal policy for ($\mathcal{P}^{CS}$) (($\mathcal{P}^{D, CS}$)) (that is, there can exist multiple  representations of $\underline\gamma^{CS, *}$ and $\underline{\tilde{\beta}}$ with $\gamma^{CS*}_{i}(y^{CS}_{i})=\tilde{\beta}^{i*}(y^{CS}_{i})$ $P$-a.s.).

 \end{itemize}
\end{prop}
\begin{proof}
Examples \ref{ex:st2} and \ref{ex:st2inv} imply Part (i), and Theorem \ref{the:4.2}  and Theorem \ref{the:csstationary policies}(i) imply Part (ii) and Part (iii), respectively. 
\end{proof}

In the following, we first establish results on the connections between uniqueness of pbp optimal policies for ($\mathcal{P}^{S}$) and ($\mathcal{P}^{D}$), which is useful in particular for LQG models, and then we applly the result to the LQG models with a partially nested information structure. The following result is a corollary to Theorems \ref{the:4.2}, \ref{the:stationary policies}, and \ref{lem:2}.

\begin{corollary}\label{corollary:un1}
Consider a stochastic dynamic team ($\mathcal{P}^{D}$) with partially nested information structure. Assume that for all $i \in \mathcal{N}$, $g_{i}$ is linear in $u^{\downarrow i}$ for all $\zeta$ (hence, Assumption \ref{assump:inv} holds). Let Assumption \ref{assump:2} hold, and let ${\underline \gamma}^{S*}\in \Gamma^{S}$ be the unique pbp optimal policy for ($\mathcal{P}^{S}$) (hence, globally optimal). Then:
\begin{itemize}[wide]
 \item [(i)] If ${\underline \gamma}^{D*}\in \Gamma^{D}$ satisfying \eqref{eq:orpds}, is affine, then ${\underline \gamma}^{D*}$ is an essentially unique affine pbp optimal policy for ($\mathcal{P}^{D}$) (unique in the class of affine policies). Moreover, if ${\hat{\underline \gamma}}^{D}\in \Gamma^{D}$ is any nonlinear stationary (pbp optimal) policy for ($\mathcal{P}^{D}$) (if it exists), then $J({{\underline \gamma}}^{D*}) \leq J({\hat{\underline \gamma}}^{D})$.
 \item [(ii)] If there exists an affine policy $\beta^{*}$ for ($\mathcal{P}^{CS}$) with representation $\beta^{*}_{i}(y^{CS}_{i}) = \gamma^{S*}_{i}(y_{i}^{S})$ for $i \in \mathcal{N}$ $P$-a.s., then $\beta^{*}$ is an essentially unique affine pbp optimal policy for ($\mathcal{P}^{CS}$) (there might exist other affine representations of the policy). Moreover, if $\hat{\beta}$ is any nonlinear pbp optimal policy for ($\mathcal{P}^{CS}$) (if it exists), then $J(\beta^{*}) \leq J(\hat{\beta})$.
 \end{itemize}
\end{corollary}

Now, we use Corollaries \ref{corollary:un} and \ref{corollary:un1} to revisit a well-known result for LQG teams (\cite{HoChu}).

\begin{example}
Consider an LQG dynamic team with  partially nested information structure where observations of DMs are given by 
\begin{flalign*}
&y_{i}^{D} = \bigg\{y^{D}_{\downarrow i}, \hat{y}^{D}_{i} :={H}_{i}\zeta + \sum_{j \in \downarrow i} B_{ij} u^{D}_{j}\bigg\},
\end{flalign*}
where $\zeta$ denotes all relevant random variables which have Gaussian distributions and mean-zero with positive covariances, and $H_{i}$ and $B_{ij}$ are matrices of appropriate dimensions. Let $I^{D}_{i}=\{{y}^{D}_{i}\}$ and the observations under the policy-dependent static reduction be given by 
$y_{i}^{S} = \{y^{S}_{\downarrow i}, \hat{y}^{S}_{i} :=H_{i}\zeta\}$. Let the expected cost function under $\underline\gamma^{D}$ be given by 
\begin{flalign*}
E^{\underline{\gamma}^{D}}[c(\omega_{0}, u_{1}, \dots, u_{N})] := E^{\underline{\gamma}^{D}}[\zeta^{\prime} Q \zeta + \underline{u}^{\prime} R \underline{u}],
\end{flalign*}
where $Q\geq 0$, $R>0$, $\underline{u}:=\{u_{1}, \dots, u_{N}\}$ and $a^{\prime}$ denotes the transpose of $a$ for $a=\zeta, \underline{u}$. Following from \cite{HoChu}, under the policy-dependent static reduction, the globally optimal policy is unique and linear (since the cost function is strictly convex in actions). Denote this globally optimal policy by  ${\underline \gamma}^{S*}=(G_{1}^{*}, \dots, G_{N}^{*})\in \Gamma^{S}$. Hence, by Proposition \ref{corollary:un}, ${\underline \gamma}^{D*}:=(K_{1}^{*}, \dots, K_{N}^{*})\in \Gamma^{D}$ with $G_{i}^{*}y_{i}^{S} = K_{i}^{*}y_{i}^{D}$ for all $i\in \mathcal{N}$ $P$-a.s. (satisfying \eqref{eq:orpds}) is the unique globally optimal policy for the dynamic team, which satisfies for all $i \in \mathcal{N}$, and $K_{i}^{*}:=(\{K_{i}^{j}\}_{j\in \downarrow i},K_{i}^{i*})$ and $G_{i}^{*}:=(\{G_{i}^{j}\}_{j\in \downarrow i},G_{i}^{i*})$,
\begin{flalign*}
&K_{i}^{j*}=G_{i}^{j*}- G_{i}^{i*}B_{ij}K_{j}^{*}~~\text{for all}~~ j \in \downarrow i,\\
&K_{i}^{i*}=G_{i}^{i*}.
\end{flalign*}

Moreover, following from Corollary \ref{corollary:un1}(i), ${\underline \gamma}^{D, *}$ is an essentially unique linear pbp optimal policy for the dynamic team.\qedwhite
\end{example}

\subsection{Convexity of Dynamic Team Problems and their Static Measurements with Control-Sharing Reduction}\label{sec:cvonx}

Next, we study convexity in policies\footnote{see, \cite[Definition 3.1]{YukselSaldiSICON17}) for the definition of convexity of a team problem in policies.} of team problems ($\mathcal{P}^{D}$), ($\mathcal{P}^{S}$), ($\mathcal{P}^{D, CS}$), and ($\mathcal{P}^{CS}$). We first present non-convex dynamic team problems (in policies) with a partially nested information structure, where the cost function satisfies Assumption \ref{assump:2}. This then shows inadequacy of Assumption \ref{assump:2}(ii) (convexity of the cost function in actions) to imply convexity of the team problem in policies for dynamic teams even for those with a partially nested information structure.

{\begin{example}\label{ex:non-cvx}
Consider a $2$-DM stochastic dynamic team ($\mathcal{P}^{D}$) where the expected cost function is given by 
\begin{flalign*}
E[c(\omega_{0}, u^{1}, u^{2})]:= E[(u^{1}+\omega_{0})^{2} + (u^{2})^{2}],
\end{flalign*}
where $\omega_{0}$ is a primitive random variable. We first discuss the convexity of the above team problem under information structure $I^{D}$ and then under $I^{CS}$ and $I^{D,CS}$.
\begin{itemize}[wide]
\item Let  $I_{1}^{D}=\{y^{D}_{1}\}$ and $I_{2}^{D}=\{y^{D}_{2}\}:=\{y^{D}_{1}, \hat{y}^{D}_{2}\}$, where $\hat{y}^{D}_{2}:=\hat{y}^{S}_{2}+u^{1}$, and $y^{D}_{1}$ and $\hat{y}^{S}_{2}$ are primitive random variables. If $\gamma^{D}_{2}(y^{D}_{2})= \sqrt[4]{\hat{y}^{D}_{2}}$, then ($\mathcal{P}^{D}$) is not convex in $(\gamma^{D}_{1}, \gamma^{D}_{1})$ since for any arbitrary policies $\gamma^{D}_{1}$ and $\hat{\gamma}^{D}_{1}$ and  for any $\alpha \in [0,1]$, we have 
\begin{flalign}
&E[c(\omega_{0}, \alpha\gamma^{D}_{1}(y^{D}_{1})+(1-\alpha)\hat{\gamma}^{D}_{1}(y^{D}_{1}),  \gamma_{2}^{D}(\hat{y}^{D}_{2})]\nonumber\\
&=E[(\omega_{0}+\alpha \gamma^{D}_{1}(y^{D}_{1}) +(1-\alpha)\hat{\gamma}^{D}_{1}(y^{D}_{1}))^{2}+\sqrt{\hat{y}^{S}_{2}+\alpha\gamma^{D}_{1}(y^{D}_{1})+(1-\alpha)\hat{\gamma}^{D}_{1}(y^{D}_{1})}]\nonumber
\end{flalign}
which implies that the above dynamic team problem ($\mathcal{P}^{D}$) is not convex in policies (this can be seen, for example by considering the trivial $\sigma$-field for DM$^{1}$, $\sigma(y^{D}_{1}):=\{\emptyset, \mathcal{F}\}$, that is DM$^{1}$ applies constant policies). We note that, under $I^{D}$ above,  DM$^{2}$  has access to  $\hat{y}^{D}_{2}$ only which is affected by the convex combination of policies of DM$^{1}$, and hence, the reduction of the observations of DM$^{2}$ is affected by the convex combination of policies of DM$^{1}$, $\alpha\gamma^{D}_{1}(y^{D}_{1})+(1-\alpha)\hat{\gamma}^{D}_{1}(y^{D}_{1})$, which may lead to non-convexity under the reduction.

\item If ${I}^{CS}_{2}:=\{y_{2}^{CS}\}:=\{y^{D}_{1}, u^{1}, \hat{y}^{S}_{2}\}$, then for any $u=\gamma^{D}_{1}(y^{D}_{1})$ and $\hat{u}=\hat{\gamma}^{D}_{1}(y^{D}_{1})$, policies $\beta_{2}^{S}$ and $\hat{\beta}_{2}^{S}$ in $\Gamma^{CS}$ can be constructed satisfying $\beta_{2}^{S}(y^{D}_{1}, \hat{y}^{S}_{2})=\gamma_{2}^{CS}(y^{D}_{1}, \hat{y}^{S}_{2}, {\gamma}^{D}_{1}(y^D_{1}))$ and $\hat{\beta}_{2}^{S}(y^{D}_{1}, \hat{y}^{S}_{2})=\gamma_{2}^{CS}(y^{D}_{1}, \hat{y}^{S}_{2}, \hat{\gamma}^{D}_{1}(y^D_{1}))$ such that
\begin{flalign*}
&{E[c(\omega_{0}, \alpha u +(1-\alpha)\hat{u}, \gamma_{2}^{CS}(y^{D}_{1}, \hat{y}^{S}_{2}, \alpha u +(1-\alpha)\hat{u}))]}\nonumber\\
&=E[(\omega_{0}+\alpha u +(1-\alpha)\hat{u})^{2}+(\alpha \beta_{2}^{S}(y^{S}_{1}, \hat{y}^{S}_{2})+(1-\alpha)\hat{\beta}_{2}^{S}(y^{S}_{1}, \hat{y}^{S}_{2}))^{2}]\nonumber,
\end{flalign*}
which implies that under the static measurements with control-sharing reduction,  the team problem above is convex in policies. We note that, under $I^{CS}$, DM$^{2}$ has access to $\hat{y}^{S}_{2}$ only without considering the convex combination of policies of DM$^{1}$ which allows for the reduction to the static information structure to be independent of policies. 

\item In view of the convexity of the team problem under $I^{CS}$, we can show that the dynamic team problem above under ${I}^{D, CS}_{2}:=\{y_{2}^{D, CS}\}:=\{y^{D}_{1}, u^{1}, \hat{y}^{D}_{2}\}$ is convex in policies using the static measurements with control-sharing reduction since the static measurements with control-sharing reduction is policy-independent. We note that, under ${I}^{D, CS}$, DM$^{2}$ has access to the convex combinations term, $\alpha u + (1-\alpha) \hat{u}$, in addition to  $\hat{y}^{D}_{2}$ which allows the DM to have access to $\hat{y}^{S}_{2}$ independent of the policies of DM$^{1}$, and this leads the reduction of the problem under ${I}^{D, CS}$ to ${I}^{CS}$ be policy-independent, and hence, convexity of the team problem under ${I}^{CS}$ leads to its convexity under ${I}^{D, CS}$. \qedwhite
\end{itemize}
\end{example}}

Now, in view of Example \ref{ex:non-cvx}, we establish below a result on convexity of team problems with partially nested information structures. We note that the following result serves as a refinement and clarification of the analysis in \cite[Section 3.3.2 and 3.3.3]{YukselSaldiSICON17}, where this distinction was not made explicit (clarification in the sense that convexity is only preserved if the information structure is partially nested with control-sharing). That is, convexity strictly requires the information structure to be partially nested with control-sharing.

\begin{theorem}\label{the:convexity}
Consider dynamic team problems ($\mathcal{P}^{D}$), ($\mathcal{P}^{S}$), ($\mathcal{P}^{CS}$), and ($\mathcal{P}^{D, CS}$) with partially nested information structure, where Assumption \ref{assump:2}(a) (on convexity of the cost function) and  Assumption \ref{assump:inv} (on the invertibility condition of observations) hold. Then:
\begin{itemize}
\item [(i)] ($\mathcal{P}^{D}$) is not necessarily convex in policies.
\item [(ii)] ($\mathcal{P}^{S}$), ($\mathcal{P}^{CS}$), and ($\mathcal{P}^{D, CS}$)  are convex in policies.

\end{itemize}
\end{theorem}
\begin{proof}
Part (i) is a straightforward consequence of Example \ref{ex:non-cvx}, and for  Part (ii), the proof for ($\mathcal{P}^{S}$) is immediate. For ($\mathcal{P}^{CS}$) and ($\mathcal{P}^{D, CS}$), since DMs can have access to actions of precedent DMs,  the reduction of each DM is independent of policies of precedent DMs (Theorem \ref{lem:2}(i)), and hence, the result follows from the fact that under the policy-dependent static reduction, the expected cost function does not change.
\end{proof}

\section{Multi-stage Team Problems: Agent-wise Optimality Analysis and Reductions}\label{sec:ms}

In this section, we consider multi-stage stochastic dynamic teams. We first provide examples, where the independent-data and AG-wise (partially) nested independent reductions apply. Then, we establish connections between AG-wise and DM-wise pbp optimal policies of dynamic multi-stage teams and their static reductions, and finally, we provide a sufficient condition under which DM-wise pbp optimality implies AG-wise pbp optimality, which leads us to use the results for the single-stage problems discussed in the previous sections.

{\begin{example}\label{lem:exmple}
Consider a multi-stage stochastic dynamic team ($\mathcal{P}^{\text{Multi}}$) with $x_{t+1}^{i}=\tilde{f}_{t}^{i}(x_{0:t}^{1:N},u_{0:t}^{1:N})+w_{t}^{i}$,
where $\tilde{f}_{t}^{i}$ is a measurable function and that observations of each AG$^{i}$ at time $t$ is of the form
$y_{t}^{i}=\tilde{h}_{t}^{i}(x_{0:t}^{1:N}, u_{0:t-1}^{1:N})+v_{t}^{i}$ with random variables $v_{t}^{i}$s being independent of other exogenous random variables of dynamics and observations, and having zero-mean Gaussian density functions ${N}_{t}^{i}$ with positive-definite covariances for all $i \in \mathcal{N}$ and $t\in \mathcal{T}$. If $I_{t}^{i}={y}^{i}_{0:t}$, then an independent-data static reduction exists. This holds since we can write 
\begin{flalign*}
y_{t}^{i}=\hat{h}_{t}^{i}(x_{0}^{1:N}, w_{0:t-1}^{1:N}, v_{0:t-1}^{1:N}, u_{0:t-1}^{1:N})+v_{t}^{i},
\end{flalign*} 
 and we can define $\phi_{t}^{i}$ and $\eta_{t}^{i}$ as
\begin{flalign*}
&\phi_{t}^{i}=\frac{{N}_{t}^{i}(y_{t}^{i}-\hat{h}_{t}^{i}(x_{0}^{1:N}, w_{0:t-1}^{1:N}, v_{0:t-1}^{1:N}, u_{0:t-1}^{1:N}))}{{N}_{t}^{i}(y_{t}^{i})}, ~~ \eta_{t}^{i}={N}_{t}^{i}(y_{t}^{i})dy_{t}^{i}.
\end{flalign*}
\end{example}
\begin{example}\label{lem:exmple1}
Consider a multi-stage stochastic dynamic team ($\mathcal{P}^{\text{Multi}}$) with $x_{t+1}^{i}=\tilde{f}_{t}^{i}(\omega_{0},x_{0:t}^{i},u_{0:t}^{i})+w_{t}^{i}$ where $\tilde{f}_{t}^{i}$ is a measurable function and  $w_{t}^{i}$ has zero-mean Gaussian density function $N_{t}^{i}$ with positive-definite covariance. Let the observations of AG$^{i}$ at time $t$ be of the form $y_{t}^{i}={h}_{t}^{i}(x_{0:t}^{i}, y_{0:t-1}^{i}, v_{0:t}^{i})$ for all $i\in \mathcal{N}$, where $\sigma(y_{t}^{i})\subset \sigma(y_{t+1}^{i})$ and $(v_{t}^{i})_{t}$ are independent of disturbances of other DMs and independent of $\omega_{0}$. If $I_{t}^{i}:=\{y_{t}^{i}\}$ for all $i=1,\dots, N$ and $t=0,\dots,T-1$, then a AG-wise nested independent reduction exists.
\end{example}}

The following corollary to Theorems \ref{lem:1}(i) and \ref{the:4.2} establishes connections between AG-wise and DM-wise pbp optimal policies of dynamic multi-stage teams and those under independent-data and AG-wise (partially) nested independent reductions.

\begin{corollary}\label{corol:multi}
Consider a multi-stage stochastic dynamic team ($\mathcal{P}^{\text{Multi}}$). 
\begin{itemize}[wide]
\item [(i)] If there exists an independent-data static reduction, then, $\underline{\pmb{\gamma}}^{*}$ is an AG-wise (DM-wise) pbp optimal policy for  ($\mathcal{P}^{\text{Multi}}$) if and only if it is an AG-wise (DM-wise) pbp optimal policy under independent-data static reduction.
\item [(ii)] If there exists a AG-wise (partially) nested independent reduction, then, $\underline{\pmb{\gamma}}^{*}$ is an AG-wise pbp optimal policy for  ($\mathcal{P}^{\text{Multi}}$) if and only if it is an AG-wise pbp optimal policy under AG-wise (partially) nested independent reduction.
\end{itemize}
\end{corollary}
We emphasize that Part (ii) is not necessarily true for DM-wise pbp optimal policies, that is, although AG-wise pbp optimal policies for (multi-stage) dynamic teams remain AG-wise pbp optimal under independent-data and AG-wise (partially) nested independent reductions, DM-wise pbp optimal policies only remain DM-wise pbp optimal  under independent-data static reductions.
\begin{proof}
Part (i) follows from Theorem \ref{lem:1}, and the fact that the independent-data static reduction is policy-independent. Part (ii) follows from the fact that in the AG-wise (partially) nested independent reduction, following from Assumption \ref{assum:decop}, the team problem can be static through agents via policy-independent static reduction, and hence, every AG-wise pbp optimal policy will be AG-wise pbp optimal under the reduction (since fixing policies of other agents, an AG-wise pbp optimal policy is globally optimal for the agent through time which will be AG-wise pbp optimal under policy-independent, policy-dependent static reductions, and static measurements with control-sharing reduction). 
\end{proof}

As we discussed earlier, every AG-wise pbp optimal policy is DM-wise pbp optimal; however, the converse statement is not true in general. In the following, we use Corollary \ref{corol:multi} to establish a variational analysis for ($\mathcal{P}^{\text{Multi}}$) under which DM-wise pbp optimal policies are AG-wise pbp optimal.

\begin{corollary}\label{multi-var}
Consider a multi-stage stochastic dynamic team ($\mathcal{P}^{\text{Multi}}$). Assume that there exists an independent-data static reduction. Let $\underline{\pmb{\gamma}}^{*}$ be a (DM-wise) pbp optimal policy for ($\mathcal{P}^{\text{Multi}}$). Assume further that, for every $t\in \mathcal{T}$,
\begin{itemize}[wide]
\item [(i)] $\hat{c}$ (see \eqref{eq:ncostms}) is continuously differentiable in $\underline{\pmb{u}}=(\pmb{u}^{1}, \dots, \pmb{u}^{N})$,
\item [(ii)] for every $i\in \mathcal{N}$, $\hat{c}$ is convex in $\pmb{u}^{i}$, where policies of other agents (for AG$^{j}$s with $j \in \{1,\dots, i-1, i+1, \dots, N\}$) are fixed to be $\pmb{\gamma}^{-i*}$.
\end{itemize}

If for all $i \in \mathcal{N}$ and for all $\pmb{\gamma}^{i}\in \bf{\Gamma}^{i}$ with $E^{\pmb{\gamma}^{i}, \pmb{\gamma}^{-i*}}[\hat{c}(\cdot)]<\infty$,
\begin{flalign*}
E\bigg[&\nabla_{u^{i}_{t}} \hat{c}\bigg(\omega_{0}, {x}_{0}, \underline{\pmb{w}},\underline{\pmb{v}}, u_{t}^{i}, \pmb{\gamma}^{-i*}(\pmb{y}^{-i}),\gamma^{i*}_{0:t-1}(y^{i}_{0:t-1}),\gamma^{i*}_{t+1:T-1}(y^{i}_{t+1:T-1}) , \underline{\pmb{y}}\bigg)\bigg|_{u_{t}^{i}=\gamma^{i*}_{t}(y^{i}_{t})} \\
&\times \bigg(\gamma^{i}_{t}(y^{i}_{t}) -\gamma^{i*}_{t}(y^{i}_{t})\bigg)\bigg] <\infty\:\:\:\:\:\text{for all}~~t \in \mathcal{T},\nonumber
\end{flalign*}
then $\underline{\pmb{\gamma}}^{*}$ is AG-wise pbp optimal for ($\mathcal{P}^{\text{Multi}}$).
\end{corollary}
\begin{proof}
Following from Corollary \ref{corol:multi}, $\underline{\pmb{\gamma}}^{*}$ is also (DM-wise) pbp optimal under an independent-data static reduction. By fixing policies $\pmb{\gamma}^{-i*}$ and using convexity and regularity conditions under the reduction, similar to \cite[Theorems 3.3 and 3.4]{YukselSaldiSICON17}, we can show that  $\underline{\pmb{\gamma}}^{*}$ is a AG-wise pbp optimal policy under an independent-data static reduction (through showing that by fixing policies $\pmb{\gamma}^{-i*}$, ${\pmb{\gamma}}^{i*}$ is globally optimal for AG$^{i}$).
\end{proof}

\begin{remark}\label{rem:NR}
For dynamic teams under a AG-wise (partially) nested independent reduction, the above variational analysis might not hold in general since if $\underline{\pmb{\gamma}}^{*}$ is (DM-wise) pbp optimal, it may not be (DM-wise) pbp optimal under a AG-wise (partially) nested independent reduction. However, for dynamic teams under a AG-wise (partially) nested independent reduction, since under Assumption \ref{assum:decop} the team is static through agents under policy-independent reduction, by considering joint perturbations through times of a given agent (through considering an AG-wise stationary policy), variational inqualities (see e.g., Corollary \ref{claim:1} and \cite[Theorems 3.3 and 3.4]{YukselSaldiSICON17}) show the global optimality of AG-wise stationary (pbp optimal) policies. Hence, under a AG-wise (partially) nested independent reduction, if the cost function is convex and continuously differentiable in actions, then, a variational analysis guarantees global optimality of AG-wise stationary (pbp optimal) policies but not DM-wise stationary (pbp optimal) policies.
\end{remark}
\section{Conclusion}\label{conc}
In this paper, we have studied connections between stationary (pbp optimal, globally optimal) policies of dynamic teams and their static reductions. We have discussed these connections for dynamic teams under both policy-independent and policy-dependent static reductions. We have showed the existence of a bijection between policies under policy-independent static reductions, and have presented some negative results as well as sufficient conditions for some positive results, where connections can be established for dynamic teams under the policy-dependent static reductions. A summary of the connections has been depicted in Fig. \ref{fig:2.1} and  Fig. \ref{fig:2.2}. In addition, we have introduced a new information structure, static measurements with control-sharing reduction, which facilitates our analysis in establishing connections between optimality concepts as well as convexity (in policies) under the reduction. A summary of connections under this reduction has been depicted in Fig. \ref{fig:3}. Furthermore, we have presented results on multi-stage team problems where two reductions have been introduced in view of AG-wise optimality concept.

For general stochastic games, information structures entail significant subtleties not present in the theory of stochastic teams. Part II  of this paper addresses  these subtleties. 

\section*{Appendix}

\subsection{Proof of Theorem \ref{lem:1}}
We first recall sufficient conditions for the {\it Bayes Formula} (e.g., \cite[p.  216]{durrett2010probability}) which is used in the proof of Theorem \ref{lem:1}.

\begin{lemma}\label{lem:Bayes}
Consider a probability space $(\hat{\Omega}, \hat{\cal{F}}, \hat{\mathbb{P}})$ where $\hat{\mathbb{P}}$ is absolutely continuous with respect to some probability measure $\hat{\mathbb{Q}}$. If a $\sigma$-field $\cal{G}\subset \hat{\cal{F}}$, and a random variable $X$ is integrable ($E_{\hat{\mathbb{P}}}[|X|]<\infty$), then the Bayes formula holds, that is, $\hat{\mathbb{P}}$-a.s
\begin{flalign*}
E_{\hat{\mathbb{P}}}[X|{\cal{G}}]=\frac{E_{\hat{\mathbb{Q}}}[X\frac{d\hat{\mathbb{P}}}{d\hat{\mathbb{Q}}}|{\cal{G}}]}{E_{\hat{\mathbb{Q}}}[\frac{d\hat{\mathbb{P}}}{d\hat{\mathbb{Q}}}|{\cal{G}}]}.
\end{flalign*}
\end{lemma}
\begin{proof}[Proof of Theorem \ref{lem:1}]
Since policies do not change under the reduction, the proof of the result for globally and pbp optimal policies follows from \eqref{eq:2.4}. We therefore prove the result for stationary policies. Let $\underline{\gamma}^{*}$ be a stationary policy for ($\mathcal{P}^{D}$). In the following, we show that if $\underline{\gamma}^{*}$ satisfies \eqref{eq:extrastationary}, then it is also stationary under a policy-independent static reduction. Since $\underline{\gamma}^{*}$ is a stationary policy for ($\mathcal{P}^{D}$), using Lemma \ref{lem:Bayes}, we have $P$-a.s.,
\begin{flalign}
0&=\nabla_{u^{i}}E^{\gamma^{-i*}}_{{{P}}}[{c}(\omega_{0}, u^{1},\dots, u^{N})|y^{i}]|_{u^{i}=\gamma^{i*}(y^{i})}\nonumber\\
&=\nabla_{u^{i}} \bigg\{\frac{E^{\gamma^{-i*}}_{{\mathbb{Q}}}[\tilde{c}(\omega_{0}, u^{1},\dots, u^{N}, y^{1}, \dots, y^{N})|y^{i}]}{E^{\gamma^{-i*}}_{{\mathbb{Q}}}[\frac{d{{P}}}{d{\mathbb{Q}}}|y^{i}]}\bigg\}\bigg |_{u^{i}=\gamma^{i*}(y^{i})}\label{eq:der1},
\end{flalign}
where \eqref{eq:der1} follows from Lemma \ref{lem:Bayes}. Hence,
\begin{flalign}
& \bigg\{\frac{\bigg(\nabla_{u^{i}}E^{\gamma^{-i*}}_{{\mathbb{Q}}}[\tilde{c}(\omega_{0}, u^{1},\dots, u^{N}, y^{1}, \dots, y^{N})|y^{i}]\bigg)E^{\gamma^{-i*}}_{{\mathbb{Q}}}[\frac{d{{P}}}{d{\mathbb{Q}}}|y^{i}]}{\bigg(E^{\gamma^{-i*}}_{{\mathbb{Q}}}[\frac{d{{P}}}{d{\mathbb{Q}}}|y^{i}]\bigg)^{2}}\nonumber\\
& - \frac{E^{\gamma^{-i*}}_{{\mathbb{Q}}}[\tilde{c}(\omega_{0}, u^{1},\dots, u^{N}, y^{1}, \dots, y^{N})|y^{i}]\bigg(\nabla_{u^{i}}E^{\gamma^{-i*}}_{{\mathbb{Q}}}[\frac{d{{P}}}{d{\mathbb{Q}}}|y^{i}]\bigg)}{\bigg(E^{\gamma^{-i*}}_{{\mathbb{Q}}}[\frac{d{{P}}}{d{\mathbb{Q}}}|y^{i}]\bigg)^{2}}\bigg\}\bigg |_{u^{i}=\gamma^{i*}(y^{i})}=0\label{eq:ex2}.
\end{flalign}
Since $\underline{\gamma}^{*}$ satisfies \eqref{eq:extrastationary}, we have $P$-a.s. the second line of \eqref{eq:ex2} is equal to zero, and since $\frac{d{{P}}}{d{\mathbb{Q}}}>0$ $P$-a.s., and the first line of \eqref{eq:ex2} must equal to zero $P$-a.s., we have $P$-a.s.
\begin{flalign}\label{eq:exred}
\nabla_{u^{i}}E^{\gamma^{-i*}}_{{\mathbb{Q}}}[\tilde{c}(\omega_{0}, u^{1},\dots, u^{N}, y^{1}, \dots, y^{N})|y^{i}]|_{u^{i}=\gamma^{i*}(y^{i})}=0,
\end{flalign} 
which implies that $\underline{\gamma}^{*}$ is a stationary policy for ($\mathcal{P}^{D}$) under policy-independent static reductions. For the converse statement, suppose a policy $\underline{\gamma}^{*}$ is stationary for ($\mathcal{P}^{D}$) under a policy-independent static reduction (that is, \eqref{eq:exred} holds) and satisfies \eqref{eq:extrastationary}. Then, \eqref{eq:ex2} is equal to zero $P$-a.s. (since its two lines are zero), which implies that \eqref{eq:der1} holds, and hence, $\underline{\gamma}^{*}$ is a stationary policy for the dynamic team ($\mathcal{P}^{D}$), and this completes the proof.
\end{proof}

\subsection{Proof of Theorem \ref{the:stationary policies}}
For simplicity of our analysis, we consider $2$-DM teams; however, a similar argument used in the proof can be utilized for $N$-DM team problems as well. Consider a $2$-DM stochastic dynamic team ($\mathcal{P}^{D}$) with $I^{1}=\{y^{D}_{1}\}$ and $I^{2}=\{y_{2}^{D}\}:=\{y^{D}_{1}, \hat{y}^{D}_{2}\}$, where $\hat{y}^{D}_{2}=g(\hat{y}^{S}_{2},u^{1})$, and $y^{D}_{1}$ and $\hat{y}^{S}_{2}$ are primitive random variables.
\begin{itemize}[wide]
\item [Part (i), ``$\Rightarrow$":] We first show that if ${\underline \gamma}^{D*}$ satisfying Condition (C) is not a stationary policy for  ($\mathcal{P}^{D}$), then ${\underline \gamma}^{S*}$ is not a stationary policy for ($\mathcal{P}^{S}$). 
\end{itemize}
\begin{itemize}[wide]
\item [\it{Step 1.}] If ${\underline \gamma}^{D*}$ satisfying Condition (C) is not a stationary policy for  ($\mathcal{P}^{D}$), then there is a set $B \subseteq \Omega$ with $P(B)>0$ such that for $\hat{\omega} \in B$ (with $y_{2}^{D}(\hat{\omega})\in \mathbb{Y}^{2}$)
\begin{flalign}
&\nabla_{u^{2}}E[c(\omega_{0},\gamma_{1}^{D*}(y_{1}^{D}), u^{2})|y_{2}^{D}]_{u^{2}=\gamma_{2}^{D*}(y_{2}^{D})}\neq 0 \label{eq:sta1}
\end{flalign}
and/or
\begin{flalign}
&\nabla_{u^{1}}E[c(\omega_{0}, u^{1}, \gamma_{2}^{D*}(y^{D}_{1}, g(\hat{y}^{S}_{2},u^{1}))|y^{D}_{1}]_{u^{1}=\gamma_{1}^{D*}(y_{1}^{D})}\neq 0\label{eq:sta2},
\end{flalign}
If \eqref{eq:sta1} holds, then since \eqref{eq:orpds} holds, and since under Assumption \ref{assump:inv}, $g_i$ is invertible, we have  on the set $B$
\begin{flalign}
&\nabla_{u^{2}}E[c(\omega_{0},\gamma_{1}^{S*}(y_{1}^{S}), u^{2})|y_{2}^{S}]_{u^{2}=\gamma_{2}^{S*}(y_{2}^{S})}\neq 0\label{eq:st0}
\end{flalign}
which contradicts the assumption that ${\underline \gamma}^{S*}$ is a stationary policy for ($\mathcal{P}^{S}$). Since $\underline{\gamma}^{D*}$ satisfies Condition (C), we have for $\epsilon_{n}\in [0,1)$ close to zero and a policy $\delta^{1}\in \Gamma^{D}_{1}$
\begin{flalign}
\gamma_{2}^{D*}\bigg(y_{2, \epsilon_{n}}^{D}\bigg)&:=\gamma_{2}^{D*}\bigg(y^{D}_{1}, g(\hat{y}^{S}_{2},\gamma_{1}^{D*}(y_{1}^{D})+\epsilon_{n}\delta^{1}(y_{1}^{D}))\bigg)\nonumber\\
&= \gamma_{2}^{D*}\bigg(y^{D}_{1}, g(\hat{y}^{S}_{2},\gamma_{1}^{D*}(y_{1}^{D}))\bigg) + \epsilon_{n} \gamma_{2}^{D*} \bigg(y^{D}_{1}, g(\hat{y}^{S}_{2},\delta^{1}(y_{1}^{D}))\bigg)\label{eq:sta4},
\end{flalign}
where $y_{2, \epsilon_{n}}^{D}:=[y^{D}_{1}, g(\hat{y}^{S}_{2},\gamma_{1}^{D*}(y_{1}^{D})+\epsilon_{n}\delta^{1}(y_{1}^{D}))]$. 

\item [\it{Step 2.}]
If \eqref{eq:sta2} holds, then on the set $B$
\begin{flalign}
&{\lim\limits_{n\rightarrow \infty}\frac{1}{\epsilon_{n}}E\bigg[}{{c}\bigg(\omega_{0},\gamma_{1}^{D*}(y_{1}^{D})+\epsilon_{n}\delta^{1}(y_{1}^{D}),\gamma_{2}^{D*}(y_{2, \epsilon_{n}}^{D})\bigg)}\nonumber\\
&\:\:\:\:\:\:\:\:\:\:\:\:\:\:\:\:\:\:\:{-{c}\bigg(\omega_{0},\gamma_{1}^{D*}(y_{1}^{D}), \gamma_{2}^{D*}(y_{2}^{D})\bigg)\bigg|y_{1}^{D}\bigg] \neq 0}\label{eq:sta3},
\end{flalign}
for a non-zero $\delta^{1}(y_{1}^{D})$. Replacing \eqref{eq:sta4} in \eqref{eq:sta3}, and using Assumption \ref{assump:2}, we can see that the expression inside the conditional expectation \eqref{eq:sta3} is convex and continuously differentiable in $\epsilon_{n}$. Hence, similar to  \cite[Theorem 2 and 3]{KraMar82}, using the extended monotone convergence theorem, we can exchange the limit and the expectation. Therefore, since \eqref{eq:orpds} holds, and also since $y_{1}^{S}=y_{1}^{D}$, on the set $B$
\begin{flalign}
E\bigg[&\lim\limits_{n\rightarrow \infty}\frac{1}{\epsilon_{n}}{c}\bigg(\omega_{0},\gamma_{1}^{S*}(y_{1}^{S})+\epsilon_{n}\delta^{1}(y_{1}^{S}),\gamma_{2}^{S*}(y^{S}_{2})+ \epsilon_{n} \hat{\delta} (y^{S}_{2})\bigg)\nonumber\\
&\:\:\:\:\:\:\:\:\:\:\:\:\:{-{c}\bigg(\omega_{0},\gamma_{1}^{S*}(y_{1}^{S}), \gamma_{2}^{S*}(y_{2}^{S})\bigg)\bigg|y_{1}^{S}\bigg] \neq 0}\label{eq:sta6},
\end{flalign}
for the policy $\hat{\delta}$ with $\hat{\delta} (y^{S}_{2}):=\gamma_{2}^{D*} (y^{D}_{1}, g(\hat{y}^{S}_{2},\delta^{1}(y_{1}^{D})))$. Under Assumption \ref{assump:2}, by the chain rule of derivative and linearity of the conditional expectation, \eqref{eq:sta6} implies that on the set $B$
\begin{flalign}
&E\bigg[\lim\limits_{n\rightarrow \infty}\frac{1}{\epsilon_{n}}{c}\bigg(\omega_{0},\gamma_{1}^{S*}(y_{1}^{S})+\epsilon_{n}\delta^{1}(y_{1}^{S}),\gamma_{2}^{S*}(y^{S}_{2})\bigg)-{c}\bigg(\omega_{0},\gamma_{1}^{S*}(y_{1}^{S}), \gamma_{2}^{S*}(y_{2}^{S})\bigg)\bigg|y_{1}^{S}\bigg]\label{eq:sta8}\\
&+ E\bigg[\lim\limits_{n\rightarrow \infty}\frac{1}{\epsilon_{n}}{c}\bigg(\omega_{0},\gamma_{1}^{S*}(y_{1}^{S}),\gamma_{2}^{S*}(y^{S}_{2})+ \epsilon_{n} \hat{\delta} (y^{S}_{2})\bigg)-{c}\bigg(\omega_{0},\gamma_{1}^{S*}(y_{1}^{S}), \gamma_{2}^{S*}(y_{2}^{S})\bigg)\bigg|y_{1}^{S}\bigg]\neq 0.\nonumber
\end{flalign}
\item [\it{Step 3.}]
If the first line of \eqref{eq:sta8} is non-zero on the set $B$, then by exchanging the limit and expectation, we have on the set $B$
\begin{flalign}
&\nabla_{u^{1}}E[c(\omega_{0},u^{1}, \gamma_{2}^{S*}(y_{2}^{S}))|y_{1}^{S}]_{u^{1}=\gamma_{1}^{S*}(y_{1}^{S})} \neq 0,  \label{eq:sta9a}
\end{flalign}
which contradicts the assumption that ${\underline \gamma}^{S*}$ is a stationary policy for ($\mathcal{P}^{S}$). If the second line of \eqref{eq:sta8} is non-zero on the set $B$, then we have on the set $B$
\begin{flalign}
&\nabla_{u^{2}}E[c(\omega_{0},\gamma_{1}^{S*}(y_{1}^{S}),u^{2})|y_{2}^{S}]_{u^{2}=\gamma_{2}^{S*}(y_{2}^{S})}\neq 0 \label{eq:sta11}.
\end{flalign}
This is because, $\sigma(y^{S}_{1}) \subseteq \sigma(y^{S}_{2})$ and by the towering property of the conditional expectation, we have on the set $B$
\begin{flalign*}
&\nabla_{u^{2}}E[E[c(\omega_{0},\gamma_{1}^{S*}(y_{1}^{S}), u^{2})|y_{2}^{S}]\hat{\delta} (y^{S}_{2})|y_{1}^{S}]_{u^{2}=\gamma_{2}^{S*}(y_{2}^{S})}\neq 0.
\end{flalign*}
Hence, on the set $B$, \eqref{eq:st0} holds, which contradicts the fact that ${\underline \gamma}^{S*}$ is a stationary policy for ($\mathcal{P}^{S}$). 
\end{itemize}

\begin{itemize}[wide]
\item [``$\Leftarrow$":] For the converse statement in Part (i), we can use similar steps. First, we note that if \eqref{eq:st0} holds on a set of positive measures $B$, then \eqref{eq:sta1} holds on the set $B$, which contradicts the assumption that ${\underline \gamma}^{D*}$ is a stationary policy for ($\mathcal{P}^{D}$). Hence, $P$-a.s.,
\begin{flalign}
&\nabla_{u^{2}}E[c(\omega_{0},\gamma_{1}^{S*}(y_{1}^{S}), u^{2})|y_{2}^{S}]_{u^{2}=\gamma_{2}^{S*}(y_{2}^{S})} = 0. \label{eq:st011}
\end{flalign}
Similar to the steps above, we can show that if ${\underline \gamma}^{D*}$ is a stationary policy for ($\mathcal{P}^{D}$), then  \eqref{eq:sta8} is equal to zero $P$-a.s. Hence, this implies that either both lines of \eqref{eq:sta8} are equal to zero $P$-a.s. or none of them is equal to zero $P$-a.s. But if the first line of \eqref{eq:sta8} is not equal to zero on the set $B$, then similar to the above we can show that \eqref{eq:st0} holds on the set $B$, which contradicts \eqref{eq:st011}, and this completes the proof. 

\end{itemize}

\subsection{Proof of Theorem \ref{lem:2}}
\begin{itemize}[wide]
\item [Part (i):] This follows from Theorem \ref{the:deptoind} since the static measurements with control-sharing reduction \eqref{eq:stmr} is policy independent, and the cost function remains unchanged under the static measurements with control-sharing reduction. For the connections between stationary policies, we have $P$-a.s.,
\begin{flalign}\label{stationar1}
0& ={\nabla_{u^{i}}E\bigg[c\bigg(\omega_{0}, (\underline\gamma^{D, CS*}_{-i}(y^{D, CS}_{-i}), u^{i})\bigg)\bigg|y^{D, CS}_{i}\bigg]\bigg|_{u^{i}=\gamma^{D,  CS*}_{i}(y^{D, CS}_{i})}}\\
& {=\nabla_{u^{i}}E\bigg[c\bigg(\omega_{0},(\underline\gamma^{CS, *}_{-i}(y^{CS}_{-i}), u^{i})\bigg)\bigg|y^{D}_{i}, \gamma^{CS*}_{\downarrow i} (y^{CS}_{\downarrow i}) \bigg]\bigg|_{u^{i}=\gamma^{CS, *}_{i}(y^{CS}_{i})}}\nonumber\\
& {=\nabla_{u^{i}}E\bigg[c\bigg(\omega_{0},(\underline\gamma^{CS, *}_{-i}(y^{CS}_{-i}), u^{i})\bigg)\bigg|y^{S}_{i}, \gamma^{CS*}_{\downarrow i} (y^{CS}_{\downarrow i}) \bigg]\bigg|_{u^{i}=\gamma^{CS, *}_{i}(y^{CS}_{i})}}\nonumber\\
& {=\nabla_{u^{i}}E\bigg[c\bigg(\omega_{0},(\underline\gamma^{CS, *}_{-i}(y^{CS}_{-i}), u^{i})\bigg)\bigg|y^{CS}_{i} \bigg]\bigg|_{u^{i}=\gamma^{CS, *}_{i}(y^{CS}_{i})}}\nonumber,
\end{flalign}
where 
\begin{flalign*}
&(\underline\gamma^{D, CS*}_{-i}(y^{D, CS}_{-i}), u^{i})\\
&:=(\gamma^{D, CS*}_{1}(y^{D, CS}_{1}), \dots, \gamma^{D, CS*}_{i-1}(y^{D, CS}_{i-1}), u^{i}, \gamma^{D, CS*}_{i+1}(y^{D, CS}_{i+1, u^{i}}), \dots,\gamma^{D, CS*}_{N}(y^{D, CS}_{N, u^{i}}),\\
&(\underline\gamma^{CS, *}_{-i}(y^{CS}_{-i}), u^{i}):=\gamma^{CS, *}_{1}(y^{CS}_{1}), \dots, \gamma^{CS, *}_{i-1}(y^{CS}_{i-1}), u^{i}, \gamma^{CS*}_{i+1}(y^{CS}_{i+1, u^{i}}), \dots,\gamma^{CS, *}_{N}(y^{CS}_{N, u^{i}}).
\end{flalign*}
The second line of \eqref{stationar1} follows from the relation \eqref{eq:stmr} since the static measurements with control-sharing reduction satisfying this relation is policy-independent. The third line of \eqref{stationar1} follows from Assumption \ref{assump:inv} since there is a bijection between $y^{D}_{i}$ and $y^{S}_{i}$, and this completes the proof.

\item [Part (ii):] Let $\underline{\gamma}^{D*}$ be a pbp optimal policy for ($\mathcal{P}^{D}$), and let $\underline{\gamma}^{D, CS*}\in \Gamma^{D, CS}$ be such that for all $i\in \mathcal{N}$, $\gamma^{D*}_{i}(y^{D}_{i})=\gamma^{D, CS*}_{i}(y^{D, CS}_{i})$ for all $u^{\downarrow i}$ $P$-a.s. A representation of policy $\underline{\gamma}^{D, CS*}$ is $\underline{\gamma}^{D*}$ itself, where for every $i\in \mathcal{N}$, the extra information $u^{\downarrow i}$ has not been used. In the following, we show that $\underline{\gamma}^{D*}$ is also pbp optimal for ($\mathcal{P}^{D, CS}$). Suppose that it is not; then there is an index $i\in \mathcal{N}$ and a policy $\beta^{i}\in \Gamma^{D, CS}_{i}$ (with $(\beta^{i}, \gamma^{D*}_{-i})\in \Gamma^{D, CS})$ such that
\begin{flalign}
&E\bigg[c\bigg(\omega_{0},\underline\gamma^{D, *}_{-i}(y^{D}_{-i}), \beta^{i}(y^{D}_{i}, \gamma^{D, *}_{\downarrow i}(y^{D}_{\downarrow i})))\bigg)\bigg]<E\bigg[c\bigg(\omega_{0},\underline\gamma^{D, *}_{-i}(y^{D}_{-i}), \gamma^{D, *}_{i}(y^{D}_{i})\bigg)\bigg]\label{eq:nec},
\end{flalign}
where $\underline\gamma^{D*}_{-i}(y^{D}_{-i}):=(\gamma^{D*}_{1}(y^{D}_{1}), \dots, \gamma^{D*}_{i-1}(y^{D}_{i-1}), \gamma^{D*}_{i+1}(y^{D}_{i+1}), \dots,\gamma^{D*}_{N}(y^{D}_{N}))$. Since for a policy $(\beta^{i}, \gamma^{D*}_{-i})\in \Gamma^{D, CS}$, there exists a policy $(\hat{\gamma}^{D}_{i}, \gamma^{D, *}_{-i})\in \Gamma^{D}$ such that $u^{i}=\beta^{i}(y^{D}_{i}, \gamma^{D*}_{\downarrow i}(y^{D}_{\downarrow i}))=\hat{\gamma}^{D}_{i}(y^{D}_{i})$ $P$-a.s. We note that $\gamma^{D*}_{-i}$ remains unchanged since the construction $\gamma^{D, CS*}_{-i}$ from $\underline\gamma^{D*}$ is independent of policies and only depends on actions which remain unchanged by the construction. Hence, \eqref{eq:nec} can be written as
\begin{flalign*}
&E\bigg[c\bigg(\omega_{0},\underline\gamma^{D*}_{-i}(y^{D}_{-i}), \hat{\gamma}^{D}_{i}(y^{D}_{i})\bigg)\bigg]<E\bigg[c\bigg(\omega_{0},\underline\gamma^{D*}_{-i}(y^{D}_{-i}), \gamma^{D*}_{i}(y^{D}_{i})\bigg)\bigg],
\end{flalign*}
which contradicts the assumption that $\underline{\gamma}^{D*}$ is pbp optimal for ($\mathcal{P}^{D}$).  Similarly, we can show the connections hold for stationary policies, and the negative result follows from Example \ref{ex:st2}.

\item [Part (iii):] Let $\underline{\gamma}^{S*}$ be pbp optimal for ($\mathcal{P}^{S}$), and let a policy $\underline{\gamma}^{CS*}\in \Gamma^{CS}$ be such that for all $i\in \mathcal{N}$, $\gamma^{S*}_{i}(y^{S}_{i})=\gamma^{CS*}_{i}(y^{CS}_{i})$ $P$-a.s. A representation of policy $\underline{\gamma}^{CS*}$ is $\underline{\gamma}^{S*}$ itself, where for every $i\in \mathcal{N}$, the extra information $u^{\downarrow i}$ has not been used. Similar to Part (ii), we can show that $\underline{\gamma}^{S*}$ is also pbp optimal for ($\mathcal{P}^{CS}$), which completes the proof.\qedwhite
\end{itemize}

\subsection{Proof of Corollary \ref{corollary:un1}}
A policy ${\underline \gamma}^{D*}$ and $g_{i}$ are affine in actions, and hence, ${\underline \gamma}^{D*}$ satisfies Condition (C). Since ${\underline \gamma}^{S*}$ is a stationary policy for ($\mathcal{P}^{S}$) and ${\underline \gamma}^{D*}$ satisfies Condition (C), by Theorem \ref{the:stationary policies}, ${\underline \gamma}^{D*}$ is a stationary policy (also pbp optimal using Theorem \ref{the:stationary policies}) for ($\mathcal{P}^{D}$). If there exists another linear stationary policy ${\tilde{\underline \gamma}}^{D*}$ for ($\mathcal{P}^{D}$), then by Theorem \ref{the:stationary policies}, ${\tilde{\underline \gamma}}^{S*}$ with $\tilde{\gamma}^{S*}_{i}(y_{i}^{S}) = \tilde{\gamma}^{D*}_{i}(y_{i}^{D})$ is a stationary policy for ($\mathcal{P}^{S}$), which contradicts the uniqueness of the stationary policy for ($\mathcal{P}^{S}$).  The second part is true since ${\underline \gamma}^{D*}$ is a globally optimal policy for ($\mathcal{P}^{S}$) by Theorem \ref{the:4.2}. Part (ii) can be shown similarly using Theorem \ref{the:csstationary policies}(ii).

\bibliographystyle{plain}

\end{document}